\documentclass{svmult}
\usepackage{makeidx}         
\usepackage{graphicx}        
\usepackage{multicol}        
\usepackage[bottom]{footmisc}

\makeindex             

\renewcommand{\theequation}{\arabic{section}.\arabic{equation}}
\newtheorem{thm}{Theorem}[section]

\newtheorem{rmk}{Remark}[section]

\newtheorem{defi}{Definition}[section]

\begin{document}

\baselineskip15pt

\title*{Pathwise Random Periodic Solutions
 of Stochastic Differential Equations}
\titlerunning{Pathwise Random Periodic Solutions}
\author{Chunrong Feng$^{1,2}$, Huaizhong Zhao$^1$, Bo Zhou$^1$}
\authorrunning{C. R. Feng, H. Z. Zhao and B. Zhou}
\institute{ Department of Mathematical Sciences, Loughborough
University, LE11 3TU, UK\and
Department of Mathematics, Shanghai Jiao Tong University,
Shanghai, 200240, China\\
 \texttt{C.Feng@lboro.ac.uk}, \texttt{H.Zhao@lboro.ac.uk}, \texttt{B.Zhou@lboro.ac.uk}}

\maketitle
\newcounter{bean}
\begin{abstract}
In this paper, we study the existence of random periodic solutions for
semilinear stochastic differential equations. We identify these
as the solutions of coupled forward-backward infinite horizon stochastic
integral equations in general cases. We then use the argument of the relative compactness of Wiener-Sobolev spaces in $C^0([0, T], L^2(\Omega))$ and generalized Schauder's fixed point theorem
to prove the existence of a solution of the coupled stochastic forward-backward
infinite horizon integral equations.
The condition on $F$ is then further weakened by
applying the coupling method of forward and backward Gronwall inequalities.
The results are also valid for stationary solutions as a special case when the period
$\tau$ can be an arbitrary number.\\

\vskip25pt
\end{abstract}

 \renewcommand{\theequation}{\arabic{section}.\arabic{equation}}

\section{Introduction}

Random dynamical systems arise in modeling many phenomena in
physics, biology, climatology, economics, etc., when uncertainties
or random influences, called noises, are taken into account. The
need for studying random dynamical systems was emphasised by Ulam and
von Neumann \cite{ul-ne} in 1945. It has been pushed since 1980s
by establishing random dynamical systems generated by random maps,
stochastic ordinary
differential equations and stochastic partial differential equations, we refer
the reader to \cite{ar},\cite{ki},\cite{ku2},\cite{liu-Qi},\cite{mo-zh-zh} and
the references therein.
Periodic solutions have been a central concept in the theory of the
deterministic dynamical system for over a century starting from Poincar\'e's seminal work \cite{poincare}. They have been studied for many
important problems arising in numerous physical problems e.g. van der Pol equations (\cite{vanderpol}),
Li\'enard equations \cite{lienard}. Now after over a century, this topic is still one of the most
interesting nonlinear phenomena to study in the theory of the deterministic dynamical systems.
Periodic behaviour
arises naturally in many real world problems e.g. in biological, enviromental and economic systems.
But these problems are often subject to random perturbations or under the influence of noises.
Needless to say, for random dynamical systems, to study
the pathwise random periodic solutions is of great importance.
Zhao and Zheng \cite {zh-zheng} started to study the problem and gave a definition
of the pathwise random periodic solutions for $C^1$-cocycles. It is well-known  that in the deterministic case,
the most powerful method to prove the existence of the periodic solution is to study the fixed point of  the
Poincar\'e map. However, for random dynamical systems,
 it is very difficult, if not impossible, to define a useful Poincar\'e map and to
find its fixed point as the trajectory does not return to the same
set with certainty. In this paper, we will study the $\tau$-periodic
solutions of $\tau$-periodic stochastic differential equations in
$R^d$:
\begin{eqnarray}\label{march17a}
du(t)&=&-Au(t)\,dt+F(t,u(t))\, dt+B_0(t)dW(t), \ \ \ \ t \geq s,\\
u(s)&=&x\in R^d. \nonumber
\end{eqnarray}
Denote $\Delta:=\{(t,s)\in R^2, s\leq t\}$. This equation generates
a semi-flow $u:\Delta\times R^d\times\Omega\to R^d$ when the
solution exists uniquely. Here $W$ is a two-sided Brownian motion on a probability space ($\Omega, {\cal F}, P$). Define $\theta: (-\infty,\infty)\times\Omega\to \Omega$ by $\theta_t\omega(s)=W(t+s)-W(t)$. Therefore ($\Omega, {\cal F}, P, (\theta_t)_{t\in R}$) is a metric dynamical system. Assume $F$ and $B_0$ satisfy: \vskip5pt

\noindent {\bf Condition (P)} {\it There exists a constant $\tau>0$
such that for any $t\in R$, $u\in R^d$
\begin{eqnarray*}
F(t,u)=F(t+\tau,u),\
B_0(t)=B_0(t+\tau).
\end{eqnarray*}}

First,  we give the definition of the random periodic solution
\begin{defi}
A random periodic solution of period $\tau$ of a semi-flow $u:
\Delta\times R^d\times\Omega\to R^d$ is an ${\cal F}$- measurable
map $\varphi:(-\infty, \infty)\times \Omega\to R^d$ such that
\begin{eqnarray}
u(t+\tau, t, \varphi(t,\omega), \omega)=\varphi(t+\tau,\omega)=\varphi(t, \theta_\tau \omega),
\end{eqnarray}
for any $t\in R$ and $\omega\in \Omega$.
\end{defi}

Instead of following the traditional geometric method of establishing the Poincar\'e mapping, for the stochastic semi-flow,  we will give a new analytical method for coupled infinite horizon forward-backward
integral equations. We will prove that the solution of the coupled forward-backward infinite horizon integral equation gives a random periodic solution of period $\tau$ and vice versa if the random periodic solution is
tempered. Under certain conditions, we can solve this coupled forward-backward infinite horizon integral equations. For this, we use a generalized Schauder's fixed point theorem and relative compactness argument in Wiener-Sobolev spaces of Malliavin derivatives. The stationary solution is also obtained
as a special case when the period $\tau$ can take an arbitrary number.
This is the case when the operators $F(t,u)$ and $B_0(t)$ do not depend on
time $t$.
In deterministic and random dynamical
systems, to find the existence of stationary solutions and random periodic solutions, and to
construct local stable and unstable manifolds near a hyperbolic
stationary point is a fundamental problem
(\cite{ar},\cite{du-lu-sc1},\cite{du-lu-sc2},\cite{li-lu},\cite{lian-lu},\cite{mo-zh-zh},\cite{ruelle2}). The stationary solution for the deterministic autonomous
parabolic differential equations actually is a solution of the corresponding elliptic
equation. This statement is not true for non-autonomous parabolic partial differential equations, even for the deterministic case with nonlinear terms periodic in time.
For stochastic differential equations or stochastic partial
differential equations with autonomous or time periodic nonlinear terms, to find a stationary solution
or a random periodic solution is a more difficult and subtle problem.
In fact, in literature,
researchers usually assume there is an invariant set or a stationary
solutions or a fixed point, then prove
invariant manifolds and stability results at a point of the
invariant set (\cite{ar},\cite{du-lu-sc1},\cite{du-lu-sc2},\cite{li-lu},\cite{mo-zh-zh},\cite{ruelle2}). So to know what the invariant set is and whether or not the invariant set is a stationary solution or a
random periodic solution or has more complicated topology
is a basic problem. In fact, for the existence of
stationary solutions, results are only known in very few cases
(\cite{kloeden},\cite{kh-ma-si},\cite{mo-zh-zh},\cite{si1},\cite{si2},\cite{zh-zh}).
Even for the stationary solution case,
researchers can only construct stable stationary solutions using the
convergence of the pullback of the solution or infinite horizon backward stochastic differential equations
(e.g. \cite{kloeden},\cite{li-zh},\cite{mattingly},\cite{zh-zh}). Our result actually gives a
general method to establish bistable stationary
solutions and random periodic solutions.
For the periodic stochastic differential equations, as far as we know, this is the first paper investigating  the random periodic solution, even though it is a very natural problem. Since Theorem \ref{aug20d}
is valid
in very general situations,
we believe the coupled
infinite horizon forward-backward stochastic integral equations (\ref{sep17a}) should be useful
in investigating random periodic solutions of many kinds of stochastic differential equations
and stochastic partial differential equations.

\section{Coupled Forward-Backward Infinite Horizon Stochastic Integral Equations and Random Periodic Solutions}\label{section3}
\setcounter{equation}{0}



We consider the semilinear stochastic
differential equation (\ref{march17a}).
Denote the solution by $u(t,s,x,\omega)$.
Let
$A$ be an $d\times d$ matrix, we can also regard it as a linear operator in ${\cal L} (R^d)$. Throughout this paper, we suppose that $A$ is hyperbolic, and $T_t=e^{-At}$ is a hyperbolic linear flow induced by $-A$. 
 So $R^d$ has a direct sum decomposition 
$$R^d=E^s\oplus E^u,$$
where $$E^s=span\{v: v {\rm \  is\ a \ generalized\  eigenvector\  for\  an\  eigenvalue} \  \lambda {\rm \ with \ }Re(\lambda)>0\},$$
$$E^u=span\{v: v {\rm \  is\ a \ generalized\  eigenvector\  for\  an\  eigenvalue} \  \lambda {\rm \ with \ }Re(\lambda)<0\}.$$
Denote $\mu_m$ the real part of an eigenvalue of $A$ with the largest negative real part, and $\mu_{m+1}$ the real part of the eigenvalue with the smallest positive real part.
We also define the projections onto
each subspace by
$$P^+:R^d\rightarrow E^s, \ P^-: R^d \rightarrow E^u.$$
Let $W(t)$, $t\in R$ be an $M$-dimensional Brownian motion and the filtered Wiener space is $(\Omega,\mathcal{F},(\mathcal{F}^t)_{t\in R},  P)$. Here ${\cal F}_s^t:=\sigma(W_u-W_v, s\leq v\leq u\leq t)$ and ${\cal F}^t:=\vee_{s\leq t}{\cal F}_s^t$.  Suppose $B_0 (s)$ is an $d\times M$ matrix and is globally bounded $\sup_{-\infty<s<\infty} ||B_0(s)||<\infty$.
  The solution of the initial value problem (\ref{march17a}) is given by the following variation of constants formula:
  \begin{eqnarray}
  u(t,s,x,\omega)=T_{t-s}x+\int _s^t T_{t-r}F(r,u(r,s,x,\omega))dr+\int _s^tT_{t-r}B_0(r)dW(r).
  \end{eqnarray}

We consider a solution of the following coupled forward-backward infinite horizon stochastic integral
equation, which is a ${\cal B} (R)\otimes\mathcal{F}$-measurable map
$Y: (-\infty,\infty)\times\Omega\rightarrow R^d$ satisfying
\begin{eqnarray}\label{sep17a}
Y(t,\omega)&=&\int _{-\infty}^tT_{t-s}P^+F(s,Y(s,\omega))ds-\int
_t^{\infty} T_{t-s}P^-F(s,Y(s,\omega))ds\nonumber
\\
&&+ (\omega)\biggl [\int _{-\infty}^tT_{t-s}P^+B_0(s)\,dW(s)\biggr ]-
(\omega)\biggl [\int _t^{\infty}T_{t-s}P^-B_0(s)\,dW(s)\biggr ]
\end{eqnarray}
for all $\omega\in \Omega$, $t\in(-\infty,\infty)$. We will give the following general theorem
which identifies the solution of the equation (\ref{sep17a}) and a random periodic
solution of stochastic differential equation (\ref{march17a}). First of all, we recall the
definition of a tempered random variable (Definition 4.1.1 in \cite{ar}):
\begin{defi}
A random variable $X: \Omega \to R^d$ is called tempered with respect to the dynamical system $\theta$ if
$$
\lim\limits _{r\to \pm \infty}{1\over |r|}\log |X(\theta _r\omega)|= 0.
$$
The random variable is called tempered from above (below) if in the above limit,
the function $\log$ is replaced by $\log^+$ ($\log^-$), the positive (negative)
part of the function $\log$.
\end{defi}

\begin{thm}\label{aug20d} Assume Condition (P). If Cauchy problem (\ref{march17a}) has a unique
solution $u(t,s,x,\omega)$ and the coupled forward-backward infinite horizon stochastic integral equation (\ref{sep17a}) has
one solution $Y: (-\infty,+\infty)\times \Omega\rightarrow R^d$ such that $Y(t+\tau,\omega)=Y(t,\theta_{\tau} \omega) \ {\rm for \ any} \
t\in R$ a.s., then $Y$ is a random periodic solution of
equation (\ref{march17a}) i.e.
\begin{eqnarray}
u(t+\tau,t, Y(t,\omega),\omega)=Y(t+\tau,\omega)=Y(t,\theta_{\tau} \omega) \ \ {\rm for \ any} \ \
t\in R
\ \ \ \ a.s.
\end{eqnarray}
Conversely, if equation (\ref{march17a}) has a random periodic solution $Y: (-\infty,+\infty)\times\Omega\rightarrow R^d$
of period $\tau$ which is tempered from above for each $t$,
then $Y$ is a solution of the coupled forward-backward infinite horizon
stochastic integral equation (\ref{sep17a}).
\end{thm}
{\bf Proof:} If equation (\ref{sep17a}) has a solution $Y(t,\omega)$,
 then for any $\tilde t\geq t$, we have
\begin{eqnarray*}
Y(\tilde t, \omega)&=&\int^{t}_{-\infty}T_{\tilde t-t}T_{t-s}P^{+}F(s,Y(s,\omega))ds-\int^{\infty}_{t}T_{\tilde t-t}T_{t-s}P^{-}F(s,Y(s,\omega))ds\\
&&+
(\omega)\int^{t}_{-\infty}T_{\tilde t-t}T_{t-s}P^{+}B_{0}(s)dW(s)-(\omega)\int^{\infty}_{t}T_{\tilde t-t}T_{t-s}P^{-}B_{0}
(s)dW(s)\\
&&
+\int^{\tilde t}_{t}T_{\tilde t-s}P^{+}F(s,Y(s,\omega))ds+\int_{t}^{\tilde t}T_{\tilde t-s}P^{-}F(s,Y(s,\omega))ds\\
&&+
(\omega)\int^{\tilde t}_{t}T_{\tilde t-s}P^{+}B_{0}(s)dW(s)+(\omega)\int_{t}^{\tilde t}T_{\tilde t-s}P^{-}B_{0}
(s)dW(s)\\
&=&
T_{\tilde t-t}Y(t,\omega)+\int_t^{\tilde t}T_{\tilde t-s}F(s,Y(s,\omega))ds
+(\omega)\int_t^{\tilde t}T_{\tilde t-s}B_0(s)dW(s).
\end{eqnarray*}
Therefore, $Y(\tilde t,\omega)$ is a solution of
(\ref{march17a}) with starting point $x=Y(t,\omega)$. Then by
the uniqueness of the solution of the initial value problem,
$$u(\tilde t,t, Y(t,\omega),\omega)=Y(\tilde t,\omega).$$
In particular, when $\tilde t=t+\tau$, and from the assumption $Y(t+\tau,\omega)=Y(t,\theta_{\tau} \omega) \ {\rm for \ any} \ t\in R$, we have 
$$u(t+\tau,t, Y(t,\omega),\omega)=Y(t+\tau,\omega)=Y(t,\theta_{\tau} \omega)$$
for all $t\in R$ and $\omega \in \Omega$.

Conversely, assume equation (\ref{march17a}) has a random periodic solution which is
also tempered from above. First note for any integer $m$,
\begin{eqnarray}
 Y(t,\omega)&=&u(t\pm m\tau,t,Y(t,\theta _{\mp m\tau}\omega), \theta _{\mp m\tau}\omega)\nonumber\\
 &=&T_{\pm m\tau}Y(t,\theta _{\mp m\tau}\omega)+\int _t^{t\pm m\tau}T_{t\pm m\tau-r}
 F(r,u(r,t,Y(t,\theta _{\mp m\tau}\omega),\theta _{\mp m\tau}\omega))dr\nonumber\\
 &&
 +\int _t^{t\pm m\tau}
 T_{t\pm m\tau-r}B_0(r)dW(r\mp m\tau).\nonumber
\end{eqnarray}
In particular,
\begin{eqnarray}
P^+ Y(t,\omega)&=&P^+
u(t+m\tau,t,Y(t,\theta _{-m\tau}\omega), \theta _{-m\tau}\omega)\nonumber\\
 &=&T_{ m\tau}P^+Y(t,\theta _{- m\tau}\omega)+\int _t^{t+ m\tau}T_{t+ m\tau-r}
P^+F(r,u(r,t,Y(t,\theta _{-m\tau}\omega),\theta _{- m\tau}\omega))dr\nonumber\\
 &&
 +\int _t^{t+ m\tau}
 T_{t+m\tau-r}P^+B_0(r)dW(r-m\tau)\nonumber
\\
&=&T_{ m\tau}P^+Y(t,\theta _{- m\tau}\omega)+\int _t^{t+ m\tau}T_{t+ m\tau-r}
P^+F(r-m\tau, Y(r-m\tau,\omega))dr\nonumber\\
 &&
 +\int _t^{t+ m\tau}
 T_{t+m\tau-r}P^+B_0(r-m\tau)dW(r-m\tau)
 \nonumber\\
&=& T_{ m\tau}P^+Y(t,\theta _{- m\tau}\omega)+\int _{t-m\tau}^tT_{t-r}
P^+F(r, Y(r,\omega))dr
 +\int _{t-m\tau}
^t T_{t-r}P^+B_0(r)dW(r)\nonumber\\
&\to &
\int _{-\infty}^tT_{t-r}
P^+F(r, Y(r,\omega))dr
 +\int _{-\infty}
^t T_{t-r}P^+B_0(r)dW(r)
\end{eqnarray}
as $m\to\infty$. One can see that the last convergence can be made first in $L^2(dP)$, so
$$P^+ Y(t,\omega)=\int _{-\infty}^tT_{t-r}
P^+F(r, Y(r,\omega))dr
 +\int _{-\infty}
^t T_{t-r}P^+B_0(r)dW(r)$$
in $L^2(dP)$, so also a.s.
Similarly
\begin{eqnarray}
P^- Y(t,\omega)&=&P^-
u(t-m\tau,t,Y(t,\theta _{m\tau}\omega), \theta _{m\tau}\omega)\nonumber\\
&=& T_{-m\tau}P^-Y(t,\theta _{m\tau}\omega)-\int ^{t+m\tau}_tT_{t-r}
P^-F(r, Y(r,\omega))dr
 -\int ^{t+m\tau}
_t T_{t-r}P^-B_0(r)dW(r)\nonumber\\
&\to &
-\int ^{+\infty}_tT_{t-r}
P^-F(r, Y(r,\omega))dr
 -\int ^{+\infty}
_t T_{t-r}P^-B_0(r)dW(r)
\end{eqnarray}
as $m\to\infty$. So we have
$$P^- Y(t,\omega)=-\int ^{+\infty}_tT_{t-r}
P^-F(r, Y(r,\omega))dr
 -\int ^{+\infty}
_t T_{t-r}P^-B_0(r)dW(r),\ \ a.s.$$
Therefore we have proved the converse part as $Y=P^+Y+P^-Y$.
\hfill \hfill $\sharp$

\begin{rmk}
Theorem \ref{aug20d} also holds in Hilbert space $H$, but with different assumptions on $A$ and $B$ and $W$.
Assume that
$A$:  $D(A)\subset H\rightarrow H$ is a closed linear operator and $T_t=e^{-At}$ is the
strongly continuous semigroup generated by $-A$.
 Let $E$ be another separable Hilbert space and
  $W(t)$, $t\in R$ be an $E$-valued Brownian motion which is defined
  on the canonical complete filtered Wiener space $(\Omega,\mathcal{F},(\mathcal{F}^t)_{t\in R},
  P)$ and with covariance in a separable Hilbert space $K$, where $K\subset
  E$ is a Hilbert-Schmidt embedding. Here ${\cal F}_s^t:=\sigma(W_u-W_v, s\leq v\leq u\leq t)$ and ${\cal F}^t:=\vee_{s\leq t}{\cal F}_s^t$. We refer readers to Chapter 4 of Da Prato and Zabczyk \cite{da-za1} for details.
  Suppose $B_0 (s)\in L_2(K,H)$
  is a Hilbert Schmidt linear operator with $\sup_{-\infty<s<\infty} ||B_0(s)||_2<\infty$.
  Moreover, let $A$ be a
self-adjoint operator on $H$ with a discrete non-vanishing spectrum
$\{\mu_n,\hspace{0.15cm}n\geq1\}$ which is bounded below and $\{e_n\}$ be the basis for $H$ consisting of eigenvectors of $A$. We have
$Ae_n=\mu_ne_n$ for $n\geq 1$.
 Assume further that $A^{-1}$ is trace-class.
Denote $\mu_m$ the largest negative
eigenvalue of $A$, and $\mu_{m+1}$ is its smallest positive
eigenvalue. Hence, we obtain an orthogonal splitting of $H$ by two
parts. One is $H^-=span\{e_1,e_2,\cdots ,e_m\}$ corresponding to the negative eigenvalues
$\{\mu_1,\mu_2,\ldots,\mu_m\}$. The other one is  $H^+=span\{e_{m+1}, e_{m+2},\cdots\}$ corresponding to the positive eigenvalues $\{\mu_n:\hspace{0.15cm}n\geq m+1\}$. And
$H$ can be written as
$$H:=H^+\oplus H^-.$$
We also define the projections onto
each subspace by
$$P^+:H\rightarrow H^+, \ P^-:H\rightarrow H^-.$$
Since $H^-$ is finite-dimensional, then $T_t|H^-$ on $H^-$ is invertible for each
$t\geq 0$. Therefore, we set
$T_{-t}:=[T_t|H^-]^{-1}$ from $H^-\rightarrow H^-$ for each $t\geq 0$. Then everything else discussed above can work the same way. 
\end{rmk}

Before we prove the existence of  the equation (\ref{sep17a}), we would like to recall the following standard notation that we will use later. We denote $C_p^\infty(R^n)$ the set of infinitely differentiable functions $f:R^n\to R$ such that $f$ and all its partial derivatives have polynomial growth. Let ${\mathcal S}$ be the class of smooth random variables $F$ that is $F=f(W(h_1),\cdots, W(h_n))$ with $n\in N$,  $h_1,\cdots, h_n\in L^2([0,T])$ and $f\in C_p^\infty(R^n)$. The derivative operator of a smooth random variable $F$ is the stochastic process $\{{\cal D}_t F,\  t\in [0,T]\}$ defined by (c.f. \cite{nuallart})
$${\cal D}_t F=\sum_{i=1}^n {{\partial f}\over {\partial x_i}}(W(h_1),\cdots,W(h_n))h_i(t).$$
We will denote ${\cal D}^{1,2}$ the domain of ${\cal D}$ in $L^2(\Omega)$, i.e. ${\cal D}^{1,2}$ is the closure of ${\mathcal S}$ with respect to the norm
 $$||F||_{1,2}^2=E|F|^2+E||{\cal D}_tF||^2_{L^2([0,T])}.$$

 Denote $C^0([0,T],L^2(\Omega))$ the set of continuous functions $f(\cdot, \omega)$ with the norm
 $$||f||^2=\sup_{t\in [0,T]} E|f(t)|^2<\infty.$$
 It's easy to check the following revised version of relative compactness of Wiener-Sobolev space in Bally-Saussereau \cite{bally} also holds. This kind of compactness as a purely random variable
 version without including time and space variables was investigated by Da Prato, Malliavin and Nualart \cite{da-mall}  and Peszat \cite{peszat} first.
 \begin{thm}\label{B-S}
Consider a sequence $(v_n)_{n\in N}$ of $C^0([0,T],L^2(\Omega))$. Suppose that:\\
(1) $v_n(t,\cdot)\in {\cal D}^{1,2}$ and $\sup\limits_{n\in N}\sup\limits_{t\in [0, T]} ||v_n(\cdot, t)||_{{1,2}}^2 <\infty$.\\
(2) There exists a constant $C>0$ such that for any $t_1, t_2\in [0, T]$, $\sup\limits_nE|v_n(t_1)-v_n(t_2)|^2 < C|t_1-t_2|.$
%
\\
 (3) (3i) There exists a constant $C$ such that for any $0<\alpha<\beta<T$, and $h \in R$ with $|h|<\min (\alpha, T-\beta)$,
 
 $\hskip 0.3cm$  and any $t_1, t_2\in [0,T]$, $\sup\limits_n\int_{\alpha}^{\beta}E |{\cal D}_{\theta+h}v_n(t_1)-{\cal D}_\theta v_n(t_2)|^2 d\theta<C(|h|+|t_1-t_2|)$.
 
 (3ii) For any $\epsilon>0$, there exist $0<\alpha<\beta<T$ such that $\sup\limits_n \sup\limits_{t\in[0, T]}\int_{[0, T]\backslash (\alpha,\beta)} E|{\cal D}_\theta v_n(t)|^2 d\theta  <\epsilon$.\\
Then $\{v_n, n\in N\}$ is relatively compact in $C^0([0,T],L^2(\Omega))$.
 \end{thm}
{\bf Proof:} Recall the Wiener chaos expansion 
$$v_n(t,\omega)=\sum_{m=0}^\infty I_m(f_n^m(\cdot,t))(\omega),$$
where $f_n^m(\cdot,t)$  are symmetric elements of $L^2([0, T]^m)$ for each $m\geq 0$. When $m=0$, $f_n^0(t)=Ev_n(t)$, so for any $t_1, t_2\in [0, T]$,
\begin{eqnarray*}
&&\sup_n\sup_{t\in [0,T]} |f_n^0(t)|\leq \sup_n\sup_{t\in [0,T]}\sqrt {E|v_n(t)|^2}<\infty,\\
&&\sup_n|f_n^0(t_1)-f_n^0(t_2)|\leq \sup_nE |v_n(t_1)-v_n(t_2)|\leq \sup_n\sqrt {E|v_n(t_1)-v_n(t_2)|^2}\leq \sqrt{C|t_1-t_2|}.
\end{eqnarray*}
So $\{f_n^0\}_{n=1}^\infty$ is relatively compact in $C^0([0, T])$.
  For each $m\geq 1$, using the same argument as in Bally-Saussereau \cite{bally}, we conclude for each fixed t, $\{f_n^m(\cdot, t)\}_{n\in N}$ is relatively compact in $L^2([0, T]^m)$. Moreover, for each $t_1, t_2\in [0,T]$, consider
\begin{eqnarray*}
\sup_n||f_n^m(\cdot, t_1)-f_n^m(\cdot,t_2)||_{L^2([0, T]^m)}^2
\leq \sup_n\int_0^T E|{\cal D}_\theta v_n(t_1)-{\cal D}_\theta v_n(t_2)|^2d\theta
\leq C|t_1-t_2|,
\end{eqnarray*}
and
\begin{eqnarray*}
\sup_n\sup_{t\in [0, T]} ||f_n^m(\cdot,t)||_{L^2([0, T]^m)}^2\leq\sup_n\sup_{t\in [0, T]} \int_0^TE|{\cal D}_\theta v_n(t)|^2d\theta<\infty.
\end{eqnarray*}
Then by Arzela-Ascoli lemma, we know that $\{f_n^m\}_{n=1}^\infty$ is relatively compact in $C^0([0, T], L^2([0, T]^m))$. Thus we can conclude $\{v_n\}_{n=1}^\infty$ is relatively compact in $C^0([0, T], L^2(\Omega))$ using the same argument as in \cite{bally}.
\hfill \hfill $\sharp$\\

We also need the following generalized Schauder's fixed point theorem to prove our theorem. The proof is refined from the proof of Schauder's fixed point theorem. So we don't claim complete originality but include it here for completeness. Note here we don't require the subset $S$ of Banach space $H$ to be closed, but impose $T$ to be continuous from $H$ to $H$ as the fixed point may not be in $S$.
\begin{thm} (Generalized Schauder's fixed point theorem)
Let $H$ be a Banach space, S be a convex subset of $H$. Assume a map $T: H\to H$ is continuous and $T(S)\subset S$
is relatively compact in $H$. Then $T$ has a fixed point in $H$.
\end{thm}
{\bf Proof:} Because $T(S)$ is relatively compact in $H$ and $H$ is a Banach space, so for any $n\in N$, there exists finite
${1\over n}$ - net $N_n=\{y_1, y_2,\cdots, y_{r_n}\}$ such that
$$T(S)\subset \bigcup_{i=1}^{r_n}B(y_i, {1\over n}),$$
where $B(y_i, {1\over n})=\{y: ||y-y_i||_H<{1\over n}\}$, $y_i\in T(S)$, $i=1, \cdots, r_n $. Denote $E_n:=span \{N_n\}$, the finite dimensional linear subspace spanned by $N_n$.

Define a map $I_n: T(S)\to co(N_n)$ by
\begin{eqnarray}\label{2.6}
I_n(y)=\sum_{i=1}^{r_n} y_i \lambda_i(y),
\end{eqnarray}
where $co(N_n)$ is the all convex combination of the elements in $N_n$, and
\begin{eqnarray*}
\lambda_i(y)={{m_i(y)}\over {\sum\limits_{i=1}^{r_n}m_i(y)}},\ \  m_i(y)=\cases{
1-n||y-y_i||_H, {\rm \ \ \ if } \ y\in B(y_i, {1\over n}),
 \cr 0, {\rm \ \hskip 2.5cm if } \ y\notin B(y_i, {1\over n}).}
\end{eqnarray*}
It's easy to see that $m_i(y)\geq 0$, and for any $y\in T(S)$, there exists an $i_0$ ($1\leq i_0\leq r_n$) such that
$y\in B(y_{i_0}, {1\over n})$, so $m_{i_0}(y)> 0$. Therefore $\lambda_i(y)$ ($1\leq i\leq r_n$) can be defined and satisfy
\begin{eqnarray}\label{2.7}
\lambda_i(y) \geq 0\   (1\leq i_0\leq r_n),\  \ \sum\limits_{i=1} ^{r_n} \lambda_i(y)=1.
\end{eqnarray}
So $I_n$ can be defined on $T(S)$ and from (\ref{2.6}) and (\ref{2.7}) we can see $I_n(y)$ is the convex combination of the elements in $N_n$.
Hence, $I_n(y)\in co(N_n)$. Moreover, for any $y\in T(S)$,
\begin{eqnarray}\label{2.8}
||I_n(y)-y||_H&=&||\sum_{i=1}^{r_n} y_i \lambda_i(y)-\sum_{i=1}^{r_n} y\lambda_i(y)||_H\nonumber\\
&\leq& \sum_{i=1}^{r_n} ||y_i - y||_H\lambda_i(y)\nonumber\\
&=& \sum_{y\in B(y_i, {1\over n})}^{r_n} ||y_i - y||_H\lambda_i(y)+ \sum_{y\notin B(y_i, {1\over n})}^{r_n} ||y_i - y||_H\lambda_i(y)\nonumber\\
& <& {1\over n}.
\end{eqnarray}

 Note that $T(S)\subset S$, $N_n\subset T(S)$ and $S$ is convex, so $co(N_n)\subset S$. Define $T_n:= I_n \circ T$. Then $T_n: co(N_n)\to co(N_n)$.
But $co(N_n)$ is a bounded closed convex subset in $E_n$, so by the Brouwer's fixed point theorem, there exists $x_n\in co(N_n)\subset S$ such that
\begin{eqnarray}\label{2.9}
T_nx_n=x_n.
\end{eqnarray}
On the other hand, $T(S)$ is relatively compact in $H$ and $H$ is complete, so there exists a subsequence $\{x_{n_k}\}\in S$ and $x\in H$ such that
\begin{eqnarray}\label{2.10}
Tx_{n_k}\to x, \ \ as\ k\to\infty.
\end{eqnarray}
From (\ref {2.8}) and (\ref {2.9}), we have
\begin{eqnarray}\label{2.11}
||x_n-x||_H&=&||T_nx_n-x||_H\nonumber\\
&\leq & ||I_nTx_n-Tx_n||_H+||Tx_n-x||_H\nonumber\\
&<&{1\over n} + ||Tx_n-x||_H.
\end{eqnarray}
Combining (\ref {2.10}) and (\ref {2.11}), we can get $x_{n_k}\to x$, as $k\to \infty$. As $T$ is continuous and also from (\ref {2.10}), we have
$$\hskip 7cm Tx=x.\hskip 7cm\sharp$$

Now we are going to prove that equation
(\ref{sep17a}) has a solution under some conditions. So according to Theorem \ref{aug20d},
this gives the existence of the random periodic solution for the stochastic evolution equation (\ref{march17a}).

\begin{thm}\label{aug20b}
Assume above conditions on $A$ and $B_0$. Let $F:(-\infty, \infty) \times R^d\to R^d$ be a continuous map, globally
bounded and the Jacobian $\nabla F(t,\cdot)$ be globally bounded,
    and $F$ and $B_0$ also satisfy Condition (P) and there exists a constant $L_1>0$ such that $||B_0(s_1)-B_0(s_2)||^2\leq L_1 |s_1-s_2|$.
Then there exists at least one ${\cal B}(R)\otimes\mathcal{F}$-measurable map
$Y: (-\infty,+\infty)\times\Omega\rightarrow R^d$ satisfying equation (\ref{sep17a})
 and  $Y(t+\tau, \omega)=Y(t, \theta_\tau\omega)$ for
any $t\in R$ and $\omega\in \Omega$.
\end{thm}
{\bf Proof:}
Firstly, define the ${\cal B}(R)\otimes\mathcal{F}$-measurable map
$Y_{1}:(-\infty,+\infty)\times\Omega\rightarrow R^d$ by
\begin{eqnarray}
Y_{1}(t,\omega)=(\omega)\int^{t}_{-\infty}T_{t-s}P^{+}B_{0}(s)dW(s)-(\omega)\int^{\infty}_{t}T_{t-s}P^{-}B_{0}(s)dW(s).
\end{eqnarray}
Then we have
\begin{eqnarray}
Y_{1}(t,\theta_\tau\omega)&=&(\theta_\tau\omega)\int^{t}_{-\infty}T_{t-s}P^{+}B_{0}(s)dW(s)-(\theta_\tau\omega)\int^{\infty}_{t}T_{t-s}P^{-}B_{0}(s)dW(s)\nonumber\\
&=&(\omega)\int^{t+\tau}_{-\infty}T_{t+\tau-s}P^{+}B_{0}(s)dW(s)-(\omega)\int^{\infty}_{t+\tau}T_{t+\tau-s}P^{-}B_{0}(s)dW(s)\nonumber\\
&=&Y_1(t+\tau,\omega).
\end{eqnarray}
Secondly, we need to solve the equation
\begin{eqnarray}\label{zhao5}
Z(t,\omega)=\int^{t}_{-\infty}T_{t-s}P^{+}F(s, Z(s,\omega)+Y_1(s,\omega)))ds-\int^{\infty}_{t}T_{t-s}P^{-}F(s,Z(s,\omega)+Y_1(s,\omega)))ds.
\end{eqnarray}
We will do this in the following several steps.

{\it Step 1} : Define
 \begin{eqnarray*}
 C_{\tau}^0((-\infty, +\infty), L^2(\Omega)):=\{f\in C^0((-\infty, +\infty), L^2(\Omega)):
 \ {\rm for \ any} \ \
t\in (-\infty,\infty),\
f(\tau+t,\omega)=f(t,\theta_{\tau}\omega) \}.
\end{eqnarray*}
For any $z\in  C_{\tau}^0((-\infty, +\infty), L^2(\Omega))$, define
\begin{eqnarray*}
{\cal M}(z)(t,\omega)=\int^{t}_{-\infty}T_{t-s}P^{+}F(s, z(s,\omega)+Y_1(s,\omega))ds-\int^{\infty}_{t}T_{t-s}P^{-}F(s,z(s,\omega)+Y_1(s,\omega))ds.
\end{eqnarray*}
We will prove ${\cal M}$ maps $ C_{\tau}^0((-\infty, +\infty), L^2(\Omega))$ to itself.
Firstly, ${\cal M} (z)(\cdot, \omega)$ is continuous.
 For this, taking any $t_1$, $t_2 \in
(-\infty, +\infty)$ with $t_1\leq t_2$, we have
\begin{eqnarray*}
&&E |{\cal M}(z)({t_1})-
{\cal M}(z)({t_2})|^2\\
&\leq& 2E\Big [\big|
\int^{t_1}_{-\infty}T_{t_1-s}P^{+}F(s,z({s})+Y_{1}({s}))ds
-\int^{t_2}_{-\infty}T_{t_2-s}P^{+}F(s,z({s})+Y_{1}({s}))ds\big|^2\\
&&+\big|
\int^{+\infty}_{t_2}T_{t_2-s}P^{-}F(s,z({s})+Y_{1}({s}))ds-
\int^{+\infty}_{t_2}T_{t_2-s}P^{-}F(s,z({s})+Y_{1}({s}))ds
\big|^2\Big].
\end{eqnarray*}
Note that there exists a constant $C\geq 1$ such that 
\begin{eqnarray*}
|T_t P^+|\leq &C e^{-{1\over 2}\mu_{m+1}t},&\ for\ all\ t\geq 0,\\
|T_t P^-|\leq &C e^{-{1\over 2}\mu_{m}t},&\ for\ all\ t\leq 0.
\end{eqnarray*}
For the first term, we have the following estimate,
\begin{eqnarray*}
&&E\big|
\int^{t_1}_{-\infty}T_{t_1-s}P^{+}F(s,z({s})+Y_{1}({s}))ds
-\int^{t_2}_{-\infty}T_{t_2-s}P^{+}F(s,z({s})+Y_{1}({s}))ds
\big|^2 \\
&\leq&2E\big|
\int^{t_1}_{-\infty}(T_{t_1-s}P^{+}-T_{t_2-s}P^{+})F(s,z(s)+Y_{1}({s}))ds\big|^2+
2E\big|
\int^{t_2}_{t_1}T_{t_2-s}P^{+}F(s,z(s)+Y_{1}({s}))ds\big|^2
\\
&\leq& 2E\int^{t_1}_{-\infty}|T_{t_{1}-s}P^+-T_{t_{2}-s}P^+|ds\cdot\int_{-\infty}^{t_1}|T_{t_{1}-s}P^+-T_{t_{2}-s}P^+|\cdot|F(s, z(s)+Y_{1}({s}))|^2ds\\
&&+2E\int^{t_2}_{t_1}|T_{t_{2}-s}P^+|^2|F(s,z(s)+Y_{1}({s}))|^2 ds(t_2-t_1)\\
&\leq& 2 \int^{t_1}_{-\infty}|T_{t_{1}-s}P^+|\cdot|(I-T_{t_{2}-t_1})P^+|ds\cdot\int_{-\infty}^{t_1}C e^{-{1\over 2}\mu_{m+1}(t_1-s)}|F(s,z(s)+Y_{1}({s}))|^2ds\\
&&+2||F||_\infty^2(t_2-t_1)^2\\
&\leq&C||F||_\infty^2(t_2-t_1){2\over {\mu_{m+1}}}+2||F||_\infty^2(t_2-t_1)^2\\
&\leq&C^{\prime}|t_2-t_1|.
\end{eqnarray*}
Here
$$||F||^2_\infty:=\sup_{t\in (-\infty,+\infty), u\in R^d}|F(t,u)|^2,
$$
and $C'$ is a generic constant throughout the paper.
By a similar argument for the second part, we also have
\begin{eqnarray*}
&&
E\big| \int^{+\infty}_{t_1}T_{t_1-s}P^{-}F(s,z(s) +
Y_{1}({s}))ds-
\int^{+\infty}_{t_2}T_{t_2-s}P^{-}F(s,z(s)+Y_{1}({s}))ds\big|^2 \\
\\
&\leq& C||F||_\infty^2(t_2-t_1)(-{2\over {\mu_{m}}})+2||F||_\infty^2(t_2-t_1)^2\\
&\leq& C'|t_2-t_1|.
\end{eqnarray*}
Therefore, by combining the two parts, we have
\begin{eqnarray*}
E| {\cal M}(z)({t_1})-
{\cal M}(z)({t_2})|^2\leq C'|t_2-t_1|.
\end{eqnarray*}
Secondly,
\begin{eqnarray*}
E|{\cal M}(z)(t)|^2&\leq & 2E\Big|\int_{-\infty}^t T_{t-s}P^{+}F(s, z(s)+Y_1(s))ds\Big|^2+2E\Big|\int^{+\infty}_t T_{t-s}P^{-}F(s, z(s)+Y_1(s))ds\Big|^2\\
&\leq&2C^2|| F||^2_{\infty}\bigg[\Big( \int^t_{-\infty}e^{-{1\over 2}(t-s)\mu_{m+1}}ds\Big)^2+\Big(\int^{+\infty}_t e^{-{1\over 2}(t-s)\mu_m} ds\Big)^2\bigg ]\\
&\leq&8C^2||F||^2_{\infty}(\frac{1}{\mu^2_{m+1}}+\frac{1}{\mu^2_{m}}).
\end{eqnarray*}
So
\begin{eqnarray*}
||{\cal M}(z)(t)||^2=\sup_{t\in (-\infty, +\infty)}E|{\cal M}(z)(t)|^2\leq 8C^2||F||^2_{\infty}(\frac{1}{\mu^2_{m+1}}+\frac{1}{\mu^2_{m}})<\infty.
\end{eqnarray*}
Thirdly,
\begin{eqnarray*}
&&
{\cal M}(z)(t,\theta _{\tau}\omega)\\
&=&\int^{t}_{-\infty}
T_{t-s}P^{+}F(s+\tau,z(s+\tau,\omega)+Y_{1}(s,\theta_{\tau}\omega))ds-\int^{+\infty}_{t}T_{t-s}P^{-}F(s+\tau,z(s+\tau,\omega)+Y_{1}
(s,\theta_{\tau}\omega))ds\\
&=&\int^{t+\tau}_{-\infty}T_{t+\tau-s}P^{+}F(s,z(s,\omega)+Y_{1}(s,\omega))ds-\int^{+\infty}_{t+\tau}T_{t+\tau-s}P^{-}F(s,z(s,\omega)+Y_{1}({s},
\omega))ds\\
&=& {\cal M}(z)(t+\tau,\omega).
\end{eqnarray*}
Therefore, we can see ${\cal M}$ maps $C_{\tau}^0((-\infty, +\infty), L^2(\Omega))$ into itself.

{\it Step 2}: To
see the continuity of the map ${\cal M}: C_{\tau}^0((-\infty, +\infty), L^2(\Omega))\to C_{\tau}^0((-\infty, +\infty), L^2(\Omega))$, consider $z_1,z_2\in
C_{\tau}^0((-\infty, +\infty), L^2(\Omega))$,
\begin{eqnarray}
&&\sup_{t\in (-\infty,+\infty)}E
|{\cal M}(z_1)(t)-{\cal M}(z_2)(t)|^2\nonumber\\
&\leq &2||\nabla F||^2_\infty\sup_{t\in
(-\infty,+\infty)}E\int^{t}_{-\infty}|T_{t-s}P^{+}|ds
\int^{t}_{-\infty}|T_{t-s}P^{+}|\cdot
|z_1(s)-z_2(s)|^2ds\nonumber
\\
&&+2||\nabla F||^2_\infty\sup_{t\in
(-\infty,+\infty)}E\int_{t}^{+\infty}|T_{t-s}P^{-}|ds
\int_{t}^{+\infty}|T_{t-s}P^{-}|\cdot
|z_1(s)-z_2(s)|^2ds\nonumber\\
&\leq & 8C^2||\nabla F||^2_\infty(\frac{1}{\mu_{m+1}^2}+\frac{1}{\mu_{m}^2})\sup_{t\in (-\infty,+\infty)}E|z_1(t)-z_2(t)|^2,
\end{eqnarray}
where
$$||\nabla F||^2_\infty:=\sup_{t\in(-\infty,\infty), u\in R^d}||\nabla F(t,u)||^2_{{\cal L}(R^d)}.$$
That is to say that ${\cal M}: C_{\tau}^0((-\infty, +\infty), L^2(\Omega))\to
C_{\tau}^0((-\infty, +\infty), L^2(\Omega))$ is a continuous map.

{\it Step 3}: Now let's define a subset of $C_{\tau}^0((-\infty,+\infty), L^2(\Omega))$ as follows:
\begin{eqnarray*}
C^0_{\tau,\alpha}((-\infty,+\infty),{\cal D}^{1,2}) :=&\{ &f\in C_{\tau}^0((-\infty,+\infty), L^2(\Omega)):\ f|_{[0,\tau)}\in C^0([0,\tau),{\cal D}^{1,2}),\\
&& i.e.\  ||f||^2=\sup_{t\in [0,\tau)}||f(t)||^2_{1,2}<\infty, {\rm \ and\  for \ any}\  t,r\in [0,\tau), \nonumber\\
&&E|{\cal D}_r f(t)|^2\leq \alpha_r(t), \sup_{s,r_1, r_2 \in [0,\tau)}{{E|{\cal D}_{r_1} f(s)-{\cal D}_{r_2} f(s)|^2}\over {|r_1-r_2|}}<\infty\}.
\end{eqnarray*}
Here $\alpha_r(t)$ is the solution of integral equation (see page 324 in \cite{polyanin}) 
\begin{eqnarray}\label{eqn2.16}
\alpha_r(t)=A\int_{r-2\tau}^{r+2\tau} e^{-\beta |t-s|}\alpha_r(s)ds +B,
\end{eqnarray}
where
\begin{eqnarray*}
&&A=10C^2||\nabla F||^2_\infty({1\over {\mu_{m+1}}}\sum_{i=0}^\infty e^{-{1\over 2}\mu_{m+1}i\tau}-{1\over {\mu_{m}}}\sum_{i=0}^\infty e^{{1\over 2}\mu_{m}i\tau}), \\
&&B=20C^2||\nabla F||^2_\infty ||B_0||^2_\infty({1\over {\mu^2_{m+1}}}+{1\over {\mu^2_{m}}}),\  \beta={1\over 2}min\{\mu_{m+1},-\mu_m\}. 
\end{eqnarray*}
This is a convex set. We will first prove that ${\cal M}$ maps $C^0_{\tau,\alpha}((-\infty,+\infty),{\cal D}^{1,2})$ into itself.
The Malliavin derivatives of $Y_1(t,\omega)$ and ${\cal M} (z)(t,\omega)$ can be calculated as:
\begin{eqnarray*}
&&\hskip -1cm {\cal D}_rY_1(t, \omega)=\left\{
\begin{array}{ll}
T_{t-r}P^+B_0(r),
 &{\rm if} \  r\leq t,\\
-T_{t-r}P^-B_0(r),   &{\rm if} \  r> t,
\end{array}
\right.
\end{eqnarray*}
and when $r\leq t$,  
\begin{eqnarray}
&& {\cal D}_r{\cal M}(z)(t, \omega)\nonumber\\
&=&\int_{-\infty}^t T_{t-s}P^+\nabla F(s, z(s,\omega)+Y_1(s,\omega)){\cal D}_r z(s,\omega)ds-\int_t^\infty T_{t-s}P^-\nabla F(s, z(s,\omega)+Y_1(s,\omega)){\cal D}_r z(s,\omega)ds\nonumber\\
&&-\int_{-\infty}^r T_{t-s}P^+\nabla F(s, z(s,\omega)+Y_1(s,\omega))T_{s-r}P^-B_0(r)ds\nonumber\\
&&+\int_r^t T_{t-s}P^+\nabla F(s, z(s,\omega)Y_1(s,\omega))T_{s-r}P^+B_0(r)ds\nonumber\\
&&-\int_t^\infty T_{t-s}P^-\nabla F(s, z(s,\omega)+Y_1(s,\omega))T_{s-r}P^+B_0(r))ds;
\end{eqnarray}
when $r> t$,
\begin{eqnarray}
&& {\cal D}_r{\cal M}(z)(t, \omega)\nonumber\\
&=&\int_{-\infty}^t T_{t-s}P^+\nabla F(s, z(s,\omega)+Y_1(s,\omega)){\cal D}_r z(s,\omega)ds-\int_t^\infty T_{t-s}P^-\nabla F(s, z(s,\omega)+Y_1(s,\omega)){\cal D}_r z(s,\omega)ds\nonumber\\
&&-\int_{-\infty}^t T_{t-s}P^+\nabla F(s, z(s,\omega)+Y_1(s,\omega))T_{s-r}P^-B_0(r)ds\nonumber\\
&&+\int^r_t T_{t-s}P^-\nabla F(s, z(s,\omega)Y_1(s,\omega))T_{s-r}P^-B_0(r)ds\nonumber\\
&&-\int_r^\infty T_{t-s}P^-\nabla F(s, z(s,\omega)+Y_1(s,\omega))T_{s-r}P^+B_0(r))ds.
\end{eqnarray}
So using Cauchy-Schwarz inequality, we have for any $z\in C^0_{\tau,\alpha}((-\infty,+\infty),{\cal D}^{1,2})$,
when $0\leq r\leq t< \tau$,
\begin{eqnarray*}
&&E|{\cal D}_r{\cal M}(z)(t)|^2\\
&\leq & 5E\int_{-\infty}^t |T_{t-s}P^+|\cdot|\nabla F(s, z(s)+Y_1(s))|^2ds\cdot\int_{-\infty}^t|T_{t-s}P^+|\cdot |{\cal D}_r z(s)|^2ds\nonumber\\ 
&&+5E\int^{\infty}_t |T_{t-s}P^-|\cdot|\nabla F(s, z(s)+Y_1(s))|^2ds\cdot\int_t^{\infty}|T_{t-s}P^-|\cdot |{\cal D}_r z(s)|^2ds\nonumber\\ 
&&+5E\int_r^t |T_{t-s}P^+|^2\cdot| \nabla F(s,z(s)+Y_1(s))|^2ds\cdot \int_r^t |T_{s-r}P^+|^2ds\cdot |B_0(r)|^2\\
&&+5E \int_{-\infty}^r |T_{t-s}P^+|^2\cdot|\nabla F(s, z(s)+Y_1(s))|^2ds\cdot\int_{-\infty}^r |T_{s-r}P^-|^2ds\cdot |B_0(r)|^2\\
&&+5E \int_t^{\infty} |T_{t-s}P^-|^2\cdot|\nabla F(s, z(s)+Y_1(s))|^2ds\cdot\int_t^{\infty} |T_{s-r}P^+|^2ds\cdot |B_0(r)|^2\\
&\leq&10C^2{1\over {\mu_{m+1}}}||\nabla F||_\infty^2\left [\int_{r-\tau}^r\sum_{i=0}^\infty e^{-{1\over 2}\mu_{m+1}(t-s+i\tau)}E|{\cal D}_r z(s)|^2ds
+\int_{r}^t e^{-{1\over 2}\mu_{m+1}(t-s)}E|{\cal D}_r z(s)|^2ds\right ]\\
&&-10C^2{1\over {\mu_{m}}}||\nabla F||_\infty^2\left [\int_t^{r+\tau} e^{-{1\over 2}\mu_{m}(t-s)}E|{\cal D}_r z(s)|^2ds+\int_{r+\tau}^{r+2\tau}\sum_{i=0}^\infty e^{-{1\over 2}\mu_{m}(t-s-i\tau)}E|{\cal D}_r z(s)|^2ds\right ]\\
&&+5C^2||\nabla F||_\infty^2\int_r^t e^{-{1\over 2}\mu_{m+1}(t-s)}ds\cdot \int_r^t e^{-{1\over 2}\mu_{m+1}(s-r)}ds\cdot |B_0(r)|^2\\
&&+5C^2||\nabla F||_\infty^2 \int_{-\infty}^re^{-{1\over 2}\mu_{m+1}(t-s)}ds\cdot\int_{-\infty}^r e^{-{1\over 2}\mu_{m}(s-r)}ds\cdot |B_0(r)|^2\\
&&+5C^2||\nabla F||_\infty^2 \int_t^{\infty} e^{-{1\over 2}\mu_{m}(t-s)}ds\cdot\int_t^{\infty} e^{-{1\over 2}\mu_{m+1}(s-r)}ds\cdot |B_0(r)|^2\\
&\leq&A\int_{r-2\tau}^{r+2\tau}e^{-\beta |t-s|}E|{\cal D}_r z(s)|^2ds+B\\
&\leq & A\int_{r-2\tau}^{r+2\tau}e^{-\beta |t-s|}\alpha_r(s)ds+B\\
&=&\alpha_r(t) .
\end{eqnarray*}
Similarly, we can prove the same result for the case when $0\leq t< r< \tau$.
Therefore, for any $r$ and $t$, we have
$$E|{\cal D}_r{\cal M}(z)(t)|^2\leq \alpha_r(t).$$
Moreover, the solution $\alpha_r(t)$ of equation (\ref{eqn2.16}) is continuous in $t$, so for $z\in C^0_{\tau,\alpha}((-\infty,+\infty),{\cal D}^{1,2})$, there exists a constant 
$M _1$ such that 
\begin{eqnarray*}
E|{\cal D}_r z(t)|^2 \leq M_1, \ and\ E|{\cal D}_r {\cal M}(z)(t)|^2 \leq M_1, \ for\  any\  t,\  r \in [0, \tau).
\end{eqnarray*}
 Suppose there exists $L_2\geq 0$ such that  for any $r_1, r_2, s\in [0,\tau)$,
\begin{eqnarray*}
{{E|{\cal D}_{r_1} z(s)-{\cal D}_{r_2} z(s)|^2}\over {|r_1-r_2|}}\leq L_2.
\end{eqnarray*}
Then we have when $0\leq r_1<r_2\leq t<\tau$,
\begin{eqnarray*}
&&{1\over {|r_1-r_2|}}E|{\cal D}_{r_1}{\cal M}(z)(t)-{\cal D}_{r_2}{\cal M}(z)(t)|^2\\
&\leq &{5\over {|r_1-r_2|}}\Big\{E\big|\int_{-\infty}^t T_{t-s}P^+\nabla F(s, z(s)+Y_1(s))({\cal D}_{r_1} z(s)-{\cal D}_{r_2} z(s))ds\big|^2\\
&&+E\big|\int^{\infty}_t T_{t-s}P^-\nabla F(s, z(s)+Y_1(s))({\cal D}_{r_1} z(s)-{\cal D}_{r_2} z(s))ds\big|^2\\
&&+E\big|\int_{r_1}^t T_{t-s}P^+\nabla F(s, z(s)+Y_1(s))T_{s-r_1}P^+B_0(r_1)ds\\
&&\hskip 2cm-\int_{r_2}^t T_{t-s}P^+\nabla F(s, z(s)+Y_1(s))T_{s-r_2}P^+B_0(r_2))ds\big|^2\\
&&+E\big|\int^{r_1}_{-\infty} T_{t-s}P^+\nabla F(s, z(s)+Y_1(s))T_{s-r_1}P^-B_0(r_1)ds\\
&&\hskip 2cm-\int_{-\infty}^{r_2} T_{t-s}P^+\nabla F(s, z(s)+Y_1(s))T_{s-r_2}P^-B_0(r_2))ds\big|^2\\
&&+E\big|\int_{t}^{\infty} T_{t-s}P^-\nabla F(s, z(s)+Y_1(s))(T_{s-r_1}P^+B_0(r_1)-T_{s-r_2}P^+B_0(r_2))ds\big|^2\Big\}\\
&:=&A_1+A_2+A_3+A_4+A_5.
\end{eqnarray*}
We will estimate them in the following, noting $D_{r_1} z(s,\omega)$, and $D_{r_2} z(s,\omega)$ are periodic in $s$, 
\begin{eqnarray*}
A_1&\leq& {5||\nabla F||_\infty^2\over {|r_1-r_2|}}\int_{-\infty}^t |T_{t-s}P^+|ds\cdot E\int_{-\infty}^t |T_{t-s}P^+|\cdot |{\cal D}_{r_1} z(s)-{\cal D}_{r_2} z(s)|^2ds\\
&\leq&{20C^2\over {\mu^2_{m+1}}}||\nabla F||_{\infty}^2L_2.
\end{eqnarray*}
Similarly,
\begin{eqnarray*}
A_2\leq{20C^2\over {\mu^2_{m}}}||\nabla F||_{\infty}^2L_2.
\end{eqnarray*}
For $A_3$,
\begin{eqnarray*}
A_3&=&{5\over {|r_1-r_2|}}E\Big|\int_{r_1}^{r_2}\nabla F(s, z(s)+Y_1(s))ds\cdot T_{t-r_1}P^+ B_0(r_1)\\
&&+\int_{r_2}^t T_{t-s}P^+\nabla F(s, z(s)+Y_1(s))(T_{s-r_1}P^+ (B_0(r_1)- B_0(r_2))\\
&&\hskip5cm +(T_{s-r_1}P^+ -T_{s-r_2}P^+ )B_0(r_2))ds\Big|^2\\
&\leq& {15\over {|r_1-r_2|}}||\nabla F||_{\infty}^2\Big [(r_2-r_1)^2||B_0||_{\infty}^2
+\int_{r_2}^t \big|T_{t-s}P^+ds\cdot \int_{r_2}^t\big|T_{s-r_1}P^+( B_0(r_1)- B_0(r_2))\big|^2ds\\
&&\hskip 3cm + \left (\int_{r_2}^t \big|T_{s-r_1}P^+-T_{s-r_2}P^+ \big |ds\big|B_0(r_2)\big|\right )^2\Big ]\\
&\leq& 15C^2||\nabla F||_{\infty}^2\Big[||B_0||_{\infty}^2\tau+\tau^2{2\over{\mu_{m+1}}}L_1+{4\over {\mu^2_{m+1}}}||B_0||_\infty^2\Big].
\end{eqnarray*}
The following estimates about $A_4$ and $A_5$ can be obtained similarly,
\begin{eqnarray*}
A_4
&\leq&15C^2||\nabla F||_{\infty}^2\Big[||B_0||_{\infty}^2\tau+(-{4\over{\mu_{m+1}\mu_m}})L_1+{4\over{\mu^2_{m}}}||B_0||_\infty^2\Big],\\
A_5
&\leq& 20C^2||\nabla F||_{\infty}^2\Big[(-{1\over{\mu_{m+1}\mu_m}})L_1+{1\over{\mu^2_{m+1}}}||B_0||_\infty^2\Big].
\end{eqnarray*}
So when $0\leq r_1<r_2\leq t <\tau$,
 \begin{eqnarray*}
{1\over {|r_1-r_2|}}E|{\cal D}_{r_1}{\cal M}(z)(t)-{\cal D}_{r_2}{\cal M}(z)(t)|^2\leq C^{\prime}.
\end{eqnarray*}
When $0\leq r_1<t<r_2<\tau$, 
\begin{eqnarray*}
&&{1\over {|r_1-r_2|}}E|{\cal D}_{r_2}{\cal M}(z)(t)-{\cal D}_{r_1}{\cal M}(z)(t)|^2\\
&\leq &{8\over {|r_1-r_2|}}\Big\{E\big|\int_{-\infty}^t T_{t-s}P^+\nabla F(s, z(s)+Y_1(s))({\cal D}_{r_2} z(s)-{\cal D}_{r_1} z(s))ds\big|^2\\
&&+E\big|\int^{\infty}_t T_{t-s}P^-\nabla F(s, z(s)+Y_1(s))({\cal D}_{r_2} z(s)-{\cal D}_{r_1} z(s))ds\big|^2\\
&&+E\big|\int_{-\infty}^{r_1} T_{t-s}P^+\nabla F(s, z(s)+Y_1(s))(T_{s-r_2}P^-B_0(r_2)-T_{s-r_1}P^-B_0(r_1))ds\big|^2\\
&&+E\big|\int^t_{r_1} T_{t-s}P^+\nabla F(s, z(s)+Y_1(s))T_{s-r_2}P^-B_0(r_2)ds\big|^2\\
&&+E\big|\int_t^{r_2} T_{t-s}P^-\nabla F(s, z(s)+Y_1(s))T_{s-r_1}P^-B_0(r_2)ds\big|^2\\
&&+E\big|\int_{r_1}^t T_{t-s}P^+\nabla F(s, z(s)+Y_1(s))T_{s-r_1}P^+B_0(r_1)ds\big|^2\\
&&+E\big|\int^{\infty}_{r_2} T_{t-s}P^-\nabla F(s, z(s)+Y_1(s))(T_{s-r_2}P^+B_0(r_2)-T_{s-r_1}P^+B_0(r_1))ds\big|^2\\
&&+E\big|\int_t^{r_2} T_{t-s}P^-\nabla F(s, z(s)+Y_1(s))T_{s-r_1}P^+B_0(r_1)ds\big|^2\Big\}\\
&\leq& 32C^2|\nabla F||_\infty^2\Big\{({1\over {\mu^2_{m+1}}}+{1\over {\mu^2_{m}}})L_2+(-{1\over{\mu_{m+1}\mu_m}}L_1+{1\over{\mu^2_{m}}}||B_0||_\infty^2)\\
&&+||B_0||_{\infty}^2\tau+(-{1\over{\mu_{m+1}\mu_m}}L_1+{1\over{\mu^2_{m+1}}}||B_0||_\infty^2)\Big\}\\
&:=&C^{\prime}.
\end{eqnarray*}
When $0\leq t\leq r_1<r_2<\tau$, similar to the case when $0\leq r_1<r_2\leq t<\tau$.
Therefore, ${\cal M}$ maps
 $C^0_{\tau,\alpha}((-\infty,+\infty),{\cal D}^{1,2})$ to itself. Now define
$${\cal M}(C^0_{\tau,\alpha}((-\infty,+\infty),{\cal D}^{1,2}))|_{[0,\tau)}:= \{f|_{[0,\tau)}: f\in {\cal M}(C^0_{\tau,\alpha}((-\infty,+\infty),{\cal D}^{1,2}))\}.$$
In order to prove  ${\cal M}(C^0_{\tau,\alpha}((-\infty,+\infty),{\cal D}^{1,2}))|_{[0,\tau)}$ is relatively compact in $C^0([0, \tau),L^2(\Omega))$, what left is to prove that 
${\cal D}_r {\cal M} (z)(t,\omega)$ is equicontinuous in $t$. We will consider several cases.
When $0\leq r<t_1<t_2<\tau$,
\begin{eqnarray*}
&&E|{\cal D}_{r}{\cal M}(z)(t_2)-{\cal D}_{r}{\cal M}(z)(t_1)|^2\\
&\leq &9\Big\{E\big|\int_{-\infty}^{t_1} (T_{t_2-s}P^+-T_{t_1-s}P^+)\nabla F(s, z(s)+Y_1(s)){\cal D}_{r} z(s)ds\big|^2\\
&&\hskip 0.5cm+E\big|\int^{t_2}_{t_1} T_{t_2-s}P^+\nabla F(s, z(s)+Y_1(s)){\cal D}_{r} z(s)ds\big|^2\\
&&+E\big|\int^{\infty}_{t_2}(T_{t_1-s}P^-- T_{t_2-s}P^-)\nabla F(s, z(s)+Y_1(s)){\cal D}_{r} z(s)ds\big|^2\\
&&\hskip 0.5cm+E\big|\int_{t_1}^{t_2} T_{t_1-s}P^-\nabla F(s, z(s)+Y_1(s)){\cal D}_{r} z(s)ds\big|^2\\
&&+E\big|\int^{t_1}_r (T_{t_1-s}P^+-T_{t_2-s}P^+)\nabla F(s, z(s)+Y_1(s))T_{s-r}P^+B_0(r)ds\big|^2\\
&&\hskip 0.5cm+E\big|\int_{t_1}^{t_2} T_{t_2-s}P^+\nabla F(s, z(s)+Y_1(s))T_{s-r}P^+B_0(r))ds\big|^2\\
&&+E\big|\int^{r}_{-\infty} (T_{t_2-s}P^+-T_{t_1-s}P^+)\nabla F(s, z(s)+Y_1(s))T_{s-r}P^-B_0(r)ds\big|^2\\
&&+E\big|\int_{t_1}^{\infty} (T_{t_1-s}P^--T_{t_2-s}P^-)\nabla F(s, z(s)+Y_1(s))T_{s-r}P^+B_0(r)ds\big|^2\\
&&\hskip 0.5cm+E\big|\int_{t_1}^{t_2} T_{t_1-s}P^-\nabla F(s, z(s)+Y_1(s))T_{s-r}P^+B_0(r)ds\big|^2\Big\}\\
&\leq&36C^2||\nabla F||_\infty^2\Big\{({1\over {\mu^2_{m+1}}}
+{1\over {\mu^2_{m}}})\alpha_1+2\alpha_1(t_2-t_1)\\
&&\hskip 0.5cm+||
B_0||^2_\infty ({1\over {\mu^2_{m+1}}}\tau+{1\over {\mu^2_{m}}}+{1\over {\mu^2_{m+1}}})+2||B_0||^2_\infty(t_2-t_1)\Big\}(t_2-t_1)\\
&\leq & C^{\prime}|t_2-t_1|,
\end{eqnarray*}
for a generic constant $C^{\prime}>0$.
When $0\leq t_1<r<t_2<\tau$ and $0\leq t_1<t_2<r<\tau$, similarly one can prove
\begin{eqnarray*}
E|{\cal D}_{r}{\cal M}(z)(t_2)-{\cal D}_{r}{\cal M}(z)(t_1)|^2
\leq C^{\prime}|t_2-t_1|.
\end{eqnarray*}
Thus,
from the above arguments, by Theorem \ref{B-S}, ${\cal M}(C^0_{\tau,\alpha}((-\infty,+\infty),{\cal D}^{1,2}))|_{[0,\tau)}$ is relatively compact in  $C^0([0,\tau),L^2(\Omega))$. Now we can prove that ${\cal M}( C_{\tau,\alpha}^0((-\infty, +\infty), {\cal D}^{1,2}))$ is relatively compact in $C_\tau^0((-\infty, +\infty), L^2(\Omega))$.
From the above, we know for any sequence ${\cal M}(z_n)\in
C_{\tau,\alpha}^0((-\infty, +\infty), {\cal D}^{1,2})$, there exists a
subsequence, still denoted by ${\cal M}(z_n)$ and $Z^*\in
C^0([0,\tau), L^2(\Omega))$ such that
\begin{eqnarray}\label{zhao1}
\sup\limits_{t\in [0,\tau)}E|{\cal M}(z_n)(t,\cdot)-Z^*(t,\cdot)|^2\to 0
\end{eqnarray}
as $n\to \infty$. Set for $\tau\leq t<2\tau$,
$$Z^*(t,\omega)=Z^*(t-\tau,\theta_{\tau}\omega).$$
Noting
\begin{eqnarray*}
{\cal M}(z_n)(t,\theta_{\tau}\omega)={\cal M}(z_n)(t+\tau,\omega),
\end{eqnarray*}
from (\ref{zhao1}), and the probability preserving property of $\theta$, we have
\begin{eqnarray*}
\sup\limits_{t\in [\tau,2\tau)}E|{\cal M}(z_n)(t,\cdot)-Z^*(t,\cdot)|^2
&=&\sup\limits_{t\in [0,\tau)}E|{\cal M}(z_n)(t+\tau ,\cdot)-Z^*(t+\tau,\cdot)|^2\\
&=&\sup\limits_{t\in [0,\tau)}E|{\cal M}(z_n)(t,\theta
_{\tau}\cdot)-Z^*(t,\theta_{\tau}\cdot)|^2\\
&=&\sup\limits_{t\in [0,\tau)}E|{\cal M}(z_n)(t,\cdot)-Z^*(t,\cdot)|^2\\
&\to& 0,
\end{eqnarray*}
Similarly one can prove that
\begin{eqnarray}
\sup\limits_{t\in [0,\tau)}E|{\cal M}(z_n)(t+m\tau ,\cdot)-Z^*(t+m\tau,\cdot)|^2=\sup\limits_{t\in [0,\tau)}E|{\cal M}(z_n)(t,\cdot)-Z^*(t,\cdot)|^2\to 0,
\end{eqnarray}
for any $m\in \{0, \pm1, \pm2,\cdots\}$. Therefore
\begin{eqnarray*}
\sup\limits_{t\in (-\infty,+\infty)}E|{\cal M}(z_n)(t,
\cdot)-Z^*(t,\cdot)|^2\to 0,
\end{eqnarray*}
as $n\to \infty$. Therefore ${\cal M}( C_{\tau,\alpha}^0((-\infty, +\infty), {\cal D}^{1,2}))$ is relatively compact in $C_\tau^0((-\infty, +\infty), L^2(\Omega))$.

{\it Stpe 4}: According to the generalized
Schauder's fixed point theorem, ${\cal M}$ has a fixed point in
$C_{\tau}^0((-\infty, +\infty), L^2(\Omega))$. That is to say
there exists a solution $Z\in
C_{\tau}^0((-\infty, +\infty), L^2(\Omega))$ of equation
(\ref{zhao5}) such that for any $t\in (-\infty,+\infty)$,
$Z(t+\tau,\omega)=Z(t,\theta_{\tau}\omega)$. Then $Y=Z+Y_1$ is the
desired solution of
(\ref{sep17a}).  Moreover, $Y(t+\tau,\omega)=Y(t,\theta_{\tau}\omega).$\hfill \hfill $\sharp$
\\

Now we consider the semilinear
stochastic differential equations with the additive noise of the form
\begin{eqnarray}\label{(bo50)}
du(t)&=&[-Au(t)+F(u(t))]dt+B_0dW(t),\\
u(0)&=&x\in R^d, \nonumber
\end{eqnarray}
for $t\geq 0$. Here $F$ and $B_0$ do not depend on time $t$, that is to say, $\tau$ in Condition (P) can be chosen as an arbitrary real number. We have a similar variation of constant representation
to (\ref{sep17a}). The difference is that for this equation, we have a cocycle $u: (-\infty,\infty)\times R^d\times \Omega\to R^d$.
Similar to Theorem \ref{aug20d}, we can prove the following theorem without giving the proof here.

\begin{thm}\label{09-27a}
Assume Cauchy problem (\ref{(bo50)}) has a unique
solution $u(t,x,\omega)$ and the coupled forward-backward infinite horizon
stochastic integral equation
\begin{eqnarray}\label{(bo1)}
Y(\omega)=\int^{0}_{-\infty}T_{-s}P^{+}F(Y(\theta_s\omega))ds-\int^{\infty}_{0}T_{-s}P^{-}F(Y(\theta_s
\omega))ds\nonumber\\
+(\omega)\int^{0}_{-\infty}T_{-s}P^{+}B_{0}dW(s)-(\omega)\int^{\infty}_{0}T_{-s}P^{-}B_{0}dW(s)
\end{eqnarray}
has one solution $Y: \Omega\rightarrow  R^d$, then $Y$ is a stationary solution of
equation (\ref{(bo50)}) i.e.
\begin{eqnarray}
u(t, Y(\omega),\omega)=Y(\theta_{t} \omega) \ \ {\rm for \ any} \ \
t\geq 0
\ \ \ \ a.s.
\end{eqnarray}
Conversely, if equation (\ref{(bo50)}) has a stationary solution $Y: \Omega\rightarrow R^d$
which is tempered from above,
then $Y$ is a solution of the coupled forward-backward infinite horizon
stochastic integral equation (\ref{(bo1)}).
\end{thm}

\begin{thm}\label{aug20a}
Assume the above conditions on $A$ and $B_0$. Let $F:  R^d\to  R^d$ be a continuous map, globally
bounded and $\nabla F$ be globally bounded.
 Then there exists at least one $\mathcal{F}$-measurable map
$Y:\Omega\rightarrow  R^d$ satisfying (\ref{(bo1)}).
\end{thm}
{\bf Proof:} Set the $\mathcal{F}$-measurable map
$Y_{1}:\Omega\rightarrow  R^d$
\begin{eqnarray}\label{(bo2)}
Y_{1}(\omega)=(\omega)\int^{0}_{-\infty}T_{-s}P^{+}B_{0}dW(s)-(\omega)\int^{\infty}_{0}T_{-s}P^{-}B_{0}dW(s).
\end{eqnarray}
Then we have
\begin{eqnarray*}
Y_{1}(\theta_t\omega)&=&(\theta_t\omega)\int^{0}_{-\infty}T_{-s}P^{+}B_{0}dW(s)-(\theta_t\omega)\int^{\infty}_{0}T_{-s}P^{-}B_{0}dW(s)\\
&=&(\omega)\int^{t}_{-\infty}T_{t-s}P^{+}B_{0}dW(s)-(\omega)\int^{\infty}_{t}T_{t-s}P^{-}B_{0}dW(s).
\end{eqnarray*}
We need to solve the equation
\begin{eqnarray}
Z(t,\omega)=\int^{t}_{-\infty}T_{t-s}P^{+}F(Z(s,\omega))+Y_1(\theta_s\omega))ds-\int^{\infty}_{t}T_{t-s}P^{-}F(Z(s,\omega)+Y_1(\theta_s\omega)))ds.
\end{eqnarray}
For this, define
\begin{eqnarray*}
C_s^0((-\infty, +\infty), L^2(\Omega))=&\{&f\in C((-\infty, \infty), L^2(\Omega)):
f(t,\omega)=f(0,\theta_t\omega) \ \ for \ all \ \ t\in(-\infty, +\infty)\}.
\end{eqnarray*}
We now define for any $z\in
C_s^0((-\infty, +\infty), L^2(\Omega))$,
\begin{eqnarray}
{\cal M}(z)(t,\omega)&=&\int^{t}_{-\infty}T_{t-s}P^{+}F(z(s,\omega)+Y_{1}(\theta_s\omega))ds\nonumber
\\
&&-\int^{+\infty}_{t}T_{t-s}P^{-}F(z(s,\omega)+Y_{1}(\theta_{s}\omega))ds.
\end{eqnarray}
It's easy to see that
\begin{eqnarray*}
{\cal M}(z)(0,\theta _t\omega)&=&\int^{0}_{-\infty}T_{-s}P^{+}F(z(s,\theta _t\omega)+Y_{1}(\theta_{s+t}\omega))ds-\int^{+\infty}_{0}T_{-s}P^{-}F(z(s,\theta
_t\omega)+Y_{1}(\theta_{s+t}\omega))ds\\
&=&\int^{0}_{-\infty}T_{-s}P^{+}F(z(s+t,\omega)+Y_{1}(\theta_{s+t}\omega))ds-\int^{+\infty}_{0}T_{-s}P^{-}F(z(s+t,\omega)+Y_{1}(\theta_{s+t}\omega))ds\\
&=&\int^{t}_{-\infty}T_{t-s}P^{+}F(z(s,\omega)+Y_{1}(\theta_{s}\omega))ds-\int^{+\infty}_{t}T_{t-s}P^{-}F(z(s,\omega)+Y_{1}(\theta_{s}\omega))ds\\
&=& {\cal M}(z)(t,\omega).
\end{eqnarray*}
By the similar method in the proof of Proposition \ref{aug20b}, we can see that  ${\cal M}$ maps $C_s^0((-\infty, +\infty), L^2(\Omega))$ to itself and ${\cal M}(\cdot)(t, \omega)$ is continuous.
We can for a fixed $T>0$ and define
\begin{eqnarray*}
 C_{T,\alpha}^0((-\infty, +\infty), {\cal D}^{1,2})&:=&\{ f\in C_s^0((-\infty,+\infty), L^{2}(\Omega)): f|_{[0,T)}\in C^0([0,T),{\cal D}^{1,2}),\\
&& i.e.\  ||f||^2=\sup_{t\in [0,T)}||f(t)||^2_{1,2}<\infty, {\rm \ and\  for \ any}\  t,r\in [0,T) \nonumber\\
&&E|{\cal D}_r f(t)|^2\leq \alpha_r (t), \sup_{s,r_1, r_2 \in [0,T)}{{E|{\cal D}_{r_1} f(s)-{\cal D}_{r_2} f(s)|^2}\over {|r_1-r_2|}}<\infty\}.
\end{eqnarray*}
Here $\alpha_r(t)$ is the solution of integral equation (see page 324 in \cite{polyanin}) 
\begin{eqnarray*}
\alpha_r(t)=A\int_{r-2T}^{r+2T} e^{-\beta |t-s|}\alpha_r(s)ds +B,
\end{eqnarray*}
where
\begin{eqnarray*}
&&A=10C^2||\nabla F||^2_\infty({1\over {\mu_{m+1}}}\sum_{i=0}^\infty e^{-{1\over 2}\mu_{m+1}iT}-{1\over {\mu_{m}}}\sum_{i=0}^\infty e^{{1\over 2}\mu_{m}iT}), \\
&&B=20C^2||\nabla F||^2_\infty ||B_0||^2_\infty({1\over {\mu^2_{m+1}}}+{1\over {\mu^2_{m}}}),\  \beta={1\over 2}min\{\mu_{m+1},-\mu_m\}. 
\end{eqnarray*}
And similar to the proof of Theorem \ref {aug20b}, we can prove that ${\cal M}$ maps $ C_{T,\alpha}^0((-\infty, +\infty), {\cal D}^{1,2})$ to itself and ${\cal M}(C_{T,\alpha}^0((-\infty, +\infty), {\cal D}^{1,2})))|_{[0,T)}$ is relatively compact in $C^0([0,T), L^2(\Omega))$.
We need to prove that ${\cal M}( C_{T,\alpha}^0((-\infty, +\infty), {\cal D}^{1,2})))$ is relatively compact in $C_s^0((-\infty, +\infty), L^2(\Omega))$.
 Note also for any sequence ${\cal M}(z_n)\in{\cal M}( C_{T,\alpha}^0((-\infty, +\infty), {\cal D}^{1,2})))$, there exists a
subsequence, still denoted by ${\cal M}(z_n)$ and $Z^*\in
C^0([0,T), L^2(\Omega))$ such that
$$E|{\cal M}(z_n)(0,\cdot)-Z^*(\cdot)|^2\to 0, \ \ {\rm as} \ \ n
\to \infty.$$
Define
\begin{eqnarray*}
Z^*(t,\omega)=Z^*(0,\theta_t\omega).
\end{eqnarray*}
Noting
\begin{eqnarray*}
{\cal M}(z_n)(0,\theta_{\tau}\omega)={\cal M}(z_n)(t,\omega),
\end{eqnarray*}
and by the probability preserving property of $\theta$, we have
\begin{eqnarray*}
\sup_{t\in (-\infty,\infty)}E|{\cal M}(z_n)(t,\cdot)-Z^*(\theta_t\cdot)|^2
&=&\sup_{t\in (-\infty,\infty)}E|{\cal M}(z_n)(0,\theta_t\cdot)-Z^*(\theta_t\cdot)|^2\\
&=&E|{\cal M}(z_n)(0,\cdot)-Z^*(\cdot)|^2\\
&\to& 0, \ \ {\rm as} \ \ n
\to \infty.
\end{eqnarray*}
 So
${\cal M}( C_{T,\alpha}^0((-\infty, +\infty), {\cal D}^{1,2})))$ is relatively compact in $C_s^0((-\infty, +\infty), L^2(\Omega))$.
Therefore, according to the generalized Schauder's fixed point theorem, ${\cal M}$ has a
fixed point in $C_s^0((-\infty, +\infty), L^2(\Omega))$. That is to
say that there exists $Z\in C_s^0((-\infty, +\infty), L^2(\Omega))$
such that for any $t\in (-\infty,+\infty)$,
$Z(t,\omega)=Z(0,\theta_t\omega)$ and
\begin{eqnarray*}
Z(0,\theta_t\omega)&=&\int^{t}_{-\infty}T_{t-s}P^{+}F(Z(0,\theta_{s}\omega)+Y_{1}(\theta_{s}\omega))ds-\int^{+\infty}_{t}T_{t-s}P^{-}F(Z(0,\theta_{s}\omega)+Y_{1}(\theta_{s}\omega))ds.
\end{eqnarray*}
Finally, we add $Y_1$ defined by the integral equation (\ref{(bo2)})
to the above equation and also assume
$$Y(\omega):=Z(0,\omega)+Y_1(\omega).$$
It's easy to see that $Y(\omega)$ satisfies (\ref{(bo1)}). \hfill $\sharp$\\

\begin{rmk} The stochastic periodic solution and stationary point may be non-unique. Note in the proof of Theorem \ref{aug20b} and Theorem \ref{aug20a}, the generalized Schauder's fixed point  argument
cannot guarantee the uniqueness of the fixed point. But in fact, the nonuniqueness is the
essence of stochastic periodic solutions or stationary solutions, since (random)
dynamical systems may have more than one random periodic solution
or stationary solution. So one should not expect in general there is only
one random periodic solution or stationary solution in a stochastic system.
Needless to say, it is interesting to study how many random periodic solutions or
stationary solutions that a stochastic differential equation can have.
\end{rmk}
\begin{rmk}
In Mohammed, Zhang and Zhao
\cite{mo-zh-zh}, they proved the existence
of the stationary solution for
a semilinear stochastic evolution equation under the conditions that $F$ satisfies the globally bounded and Lipschitz
conditions and the Lipschitz constant $L$ is with the restriction
$$L[\mu^{-1}_{m+1}-\mu^{-1}_{m}]<1.$$
The condition was imposed due to the use of the Banach fixed point theorem. It is noted that a similar condition appeared in many literature on random invariant manifolds of random dynamical systems (\cite{ar},\cite{du-lu-sc1},\cite{du-lu-sc2}). Therefore it would be interesting to investigate whether or not the method introduced in this paper can be used to get rid of the critical condition in the invariant manifold theorems.
\end{rmk}


\section{Weaken the Condition on $F$}
\setcounter{equation}{0}

Our purpose in this section is to push the results of last section further to find a
weaker condition to replace the global boundedness condition on $F$. For simplicity, we only consider the case when $A$ is symmetric.
We adopt all notations from the last section. 
Now we consider the following equation with a standard cut off function
$F_N^*$,
\begin{eqnarray}
\label{(bo4)}
z_N(t)&=&\int^{t}_{-\infty}T_{t-s}P^{+}F_N^*(s,z_N(s))ds-\int^{\infty}_{t}T_{t-s}P^{-}F_N^*(s,z_N(s))ds.
\end{eqnarray}
Here the function $F_N^*$  can be constructed in the following way. Denote for any $x\in E^s, y\in E^u$,
 $F^*(s,x,y):=F(s,x+Y_1^+(s),y+Y_1^-(s))$, and the cutting function is defined as following
\begin{eqnarray*}
F^*_N(s,x, y):=F^*(s,x{|x|\wedge N\over |x|},y{|y|\wedge N\over |y|}).
\end{eqnarray*}
It is easy to see that $F_N^*$ is a function from $ (-\infty,\infty)\times R^d\rightarrow R^d$ and $F_N^*$
is bounded no matter whether or not $F$ is bounded. By the previous proof, we have, as $F_N^*$
is bounded, there exists at least one ${z_N(t)}$, and the
solution depends on $N$, such that
$$||z_N||^2\leq \beta_N,$$
where $\beta_N$ is the radius of a closed ball which depends on $N$ and
is dominated by $F_N^*$ such that $$\beta_N:=8C^2||F^*_N||_\infty^2(\frac{1}{\mu^2_{m+1}}+\frac{1}{\mu^2_m}).$$
The idea here is that if we can prove there exists $\beta ^\prime>0$ which is independent of $N$, such that
for all $N$,
$$||z_N||^2\leq \beta'.$$
this is to say we can always choose $N$ big enough such that $$F_N^*=F^*,$$ and the globally bounded condition
for $F$ will then be possible to be omitted. In the following, we are going to work out the idea. To simplify the notation, we denote $z_N(t)$ by $z(t)$ in (\ref{(bo4)}) without any confusion. Consider that $F^*$ (therefore $F$) satisfies the following condition:
\medskip

\noindent
{\bf Condition (M)} {\it For any $x\in E^s$, $y\in E^u$, there exist constants $L_i>0, i=1,\cdots, 4$, and random variables $ A_1, B_1\geq 0$ such that
\begin{eqnarray*}
(x,P^+F^*(s,x,y))\leq L_1x^2+L_2y^2+A_1,\\
(y,(-P^-)F^*(s,x,y))\leq L_3x^2+L_4y^2+B_1.
\end{eqnarray*}
Set $$L_1^*=2L_1, \ L_2^*=L_2\vee {{L_2^2}\over {L_1}}, \ L_3^*=L_3\vee {{L_3^2}\over {L_4}}, \ L_4^*=2L_4,$$
 and $L_i^*, i=1,\cdots, 4$ satisfy
$$L_1^*<\mu_{m+1}, L_4^*<-\mu_m, \ L_2^*L_3^*\leq {1\over 2}(\mu_{m+1}-\mu_m-L_1^*-L_4^*)\cdot\min\{\mu_{m+1}-L_1^*, \mu_m+L_4^*\}.$$
}
Note, if we take $N$ sufficiently large such that ${{A_1}\over {2N^2}}\leq L_1$, we have that for any $x\in E^s$, $y\in E^u$, 
 \begin{eqnarray}
(x,P^+F_N^*(s,x+a,y+b))\leq L_1^*x^2+L_2^*y^2+A_1,\label{march1}\\
(y,(-P^-)F_N^*(s,x+a,y+b))\leq L_3^*x^2+L_4^*y^2+B_1\label{march2}.
\end{eqnarray}
\begin{thm}\label{09-27}
Assume conditions on $A$, $B_{0}$ in
Theorem \ref{aug20b}. Let $F:(-\infty, \infty) \times  R^d\to  R^d$ be a continuous map satisfying Condition (M), $\nabla F(t,u)$ exist and be locally bounded, and $B_0$ and $F$ satisfy Condition (P). Then the stochastic integral equation (\ref{sep17a}) has at least one solution $Y: (-\infty,+\infty)\times \Omega\rightarrow  R^d$ and the solution satisfies
\begin{eqnarray}
u(t+\tau,t, Y(t,\omega),\omega)=Y(t+\tau,\omega)=Y(t,\theta_{\tau} \omega) \ \ {\rm for \ any} \ \
t\geq 0
\ \ \ \ a.s.
\end{eqnarray}
\end{thm}
{\bf Proof:}
Consider integral equation (\ref{(bo4)}). For $z\in R^d$ and a given $Y_1\in R^d$, we have
\begin{eqnarray*}
z(t):=(z^+(t),z^-(t)),\
Y_1(t):=(Y^+_1(t),Y^-_1(t))
\end{eqnarray*}
where $z^+(t),Y^+(t)\in E^s$ and $z^-(t),Y^-(t)\in E^u$. Then
(\ref{(bo4)}) can be expressed by two parts
\begin{eqnarray}\label{(bo7)}
z^+(t)=\int^t_{-\infty}T_{t-s}P^+F_N^*(s,z^+(s),
z^-(s))ds
\end{eqnarray}
and
\begin{eqnarray}\label{(bo8)}
z^-(t)=-\int^\infty_{t}T_{t-s}P^-F_N^*(s,z^+(s),
z^-(s))ds.
\end{eqnarray}
Under the basis $\{e_i, 1\leq i\leq d\}$, we assume
\begin{eqnarray*}
z^-_i(t)&=&(z^-(t), e_i),\ i=1,\cdots, m; \\
z^+_{j}(t)&=&(z^+(t), e_{j}),\ j=m+1,\ m+2, \cdots,d.
\end{eqnarray*}
Consider the differential forms of (\ref{(bo7)}) and (\ref{(bo8)})
according to each eigenvalue of $A$, we have
\begin{eqnarray*}
\frac{dz^-_i(t)}{dt}&=&-\mu_iz^-_i(t)+(F_N^*(t,z^+(t),z^-(t)),e_i), \ i=1,\cdots, m,\\\\
\frac{dz^+_{j}(t)}{dt}&=&-\mu_{j}z^+_{j}(t)+(F_N^*(t,z^+(t),z^-(t)),e_{j}),\ j=m+1,\ m+2, \cdots, d.
\end{eqnarray*}
For the first $m$ differential equations, we consider the backward integral equations. For the
rest inequalities, we consider the forward integral equations. Then
we have
\begin{eqnarray*}
(z^-_i)^2(t)&\leq& \int^{\infty}_t 2\mu_m (z^-_i)^2(s)ds - \int^{\infty}_t 2z^-_i(s)(F_N^*(s,z^+(s),z^-(s)),e_i)ds,\\
(z^+_{j})^2(t)&\leq& \int_{-\infty}^t -2\mu_{m+1} (z^+_{j})^2(s)ds + \int_{-\infty}^t 2z^+_{j}(s)(F_N^*(s,z^+(s),z^-(s)),e_{j})ds,
\end{eqnarray*}
for $\ i=1,\cdots, m, \ j=m+1,\ m+2, \cdots,d$.
Then applying the Gronwall inequality for each
differential inequality, we have
\begin{eqnarray}
&&(z^-_i)^2(t)\leq-2\int^{\infty}_te^{-(t-s)2\mu_m}z^-_i(s)(F_N^*(s,z^+(s),z^-(s)),e_i)ds,\label{3.4}\\
&&(z^+_{j})^2(t)\leq2\int_{-\infty}^te^{-(t-s)2\mu_{m+1}}z^+_{j}(s)(F_N^*(s,z^+(s),z^-(s)),e_{j})ds\label{3.5},
\end{eqnarray}
where $\ i=1,\cdots, m, \ j=m+1,\ m+2, \cdots,d$.
Now we combine them into two types by writing
\begin{eqnarray*}
( z^+ )^2(t)=\sum_{j=m+1}^d(z^+_{j})^2(t),\ ( z^- )^2(t)=\sum_{i=1}^m(z^-_{i})^2(t)
\end{eqnarray*}
and
\begin{eqnarray*}
(z^+(s),P^+F_N^*(s,z^+(s), z^-(s)))
&=&\sum_{j=m+1}^d z^+_{j}(s)(F_N^*(s,z^+(s), z^-(s)),e_{j})\\
(z^-(s),P^-F_N^*(s,z^+(s),
z^-(s))
&=&\sum_{i=1}^m z^-_{i}(s)(F_N^*(s,z^+(s),
z^-(s)),e_{i})
\end{eqnarray*}
Then from (\ref{3.4}) and (\ref{3.5}), we have
\begin{eqnarray*}
( z^+)^2(t)&\leq& 2\int^t_{-\infty}e^{-(t-s)2\mu_{m+1}}(z^+(s),P^+F_N^*(s,z^+(s), z^-(s)))ds,\label{(bo9)}\\
( z^-)^2(t)&\leq&
-2\int^{\infty}_{t}e^{-(t-s)2\mu_m}(z^-(s),P^-F_N^*(s,z^+(s),
z^-(s)))ds.  \label{(bo10)}
\end{eqnarray*}
By (\ref{march1}) and ({\ref{march2}), we have
\begin{eqnarray}
( z^+)^2(t)\leq 2\int^t_{-\infty}e^{-(t-s)2\mu_{m+1}}[L^*_1(
z^+)^2(s)+L^*_2( z^-)^2(s)+A_1]ds\nonumber
\end{eqnarray}
\begin{eqnarray}
( z^-)^2(t)\leq 2\int^{\infty}_{t}e^{-(t-s)2\mu_m}[L^*_3(
z^+)^2(s)+L^*_4( z^-)^2(s)+B_1]ds.\nonumber
\end{eqnarray}
This will lead to
\begin{eqnarray}\label{(bo11)}
( z^+)^2(t)\leq 2\int^t_{-\infty}e^{-(t-s)2\mu_{m+1}}[L^*_1(
z^+)^2(s)+L^*_2( z^-)^2(s)]ds+\frac{A_1}{\mu_{m+1}}
\end{eqnarray}
\begin{eqnarray}\label{(bo12)}
( z^-)^2(t)\leq 2\int^{\infty}_{t}e^{-(t-s)2\mu_m}[L^*_3(
z^+)^2(s)+L^*_4( z^-)^2(s)]ds-\frac{B_1}{\mu_{m}}.
\end{eqnarray}
In the next step we will apply the Gronwall
inequality and coupling method. This leads to
\begin{eqnarray}
e^{t2\mu_{m+1}}( z^+)^2(t)&\leq&
\int^t_{-\infty}e^{s2\mu_{m+1}}2L^*_2( z^-)^2(s)ds+\frac{A_1}{\mu_{m+1}}e^{t2\mu_{m+1}}+\int^t_{-\infty}[e^{s2\mu_{m+1}}( z^+)^2(s)]2L^*_1ds.\nonumber
\end{eqnarray}
Then applying the Gronwall inequality to the above
inequality, we immediately have
\begin{eqnarray*}
e^{t2\mu_{m+1}}( z^+)^2(t)\leq
\int^t_{-\infty}e^{s2\mu_{m+1}}2L^*_2( z^-)^2(s)e^{2L^*_1(t-s)}ds+\int^t_{-\infty}2A_1e^{s2\mu_{m+1}}e^{2L^*_1(t-s)}ds.
\end{eqnarray*}
So it is trivial to see that
\begin{eqnarray}\label{(bo15)}
( z^+)^2(t)&\leq&
\int^t_{-\infty}e^{(t-s)2(L^*_1-\mu_{m+1})}2L^*_2( z^-)^2(s)ds+2A_1\int^t_{-\infty}e^{(t-s)2(L^*_1-\mu_{m+1})}ds\nonumber\\
&\leq & 2\int^t_{-\infty}e^{(t-s)2(L^*_1-\mu_{m+1})}L^*_2(
z^-)^2(s)ds+\frac{A_1}{\mu_{m+1}-L^*_1}.
\end{eqnarray}
From (\ref{(bo12)}) we have
\begin{eqnarray*}
e^{t2\mu_{m}}( z^-)^2(t)&\leq&
\int^{\infty}_{t}e^{s2\mu_{m}}2L^*_3( z^+)^2(s)ds-\frac{B_1}{\mu_{m}}e^{t2\mu_m}+\int^{\infty}_{t}[e^{s2\mu_{m}}( z^-)^2(s)]2L^*_4ds.
\end{eqnarray*}
Applying the Gronwall inequality, we have
\begin{eqnarray*}
e^{t2\mu_{m}}(z^-)^2(t)&\leq&
\int^{\infty}_{t}e^{s2\mu_{m}}2L^*_3( z^+)^2(s)e^{2L^*_4(s-t)}ds+\int^{\infty}_{t}2B_1e^{s2\mu_{m}}e^{2L^*_4(s-t)}ds.
\end{eqnarray*}
So it is trivial to see that
\begin{eqnarray}\label{(bo16)}
( z^-)^2(t)&\leq&
\int^{\infty}_{t}e^{(s-t)2(\mu_{m}+L^*_4)}2L^*_3( z^+)^2(s)ds+2B_1\int^{\infty}_{t}e^{(s-t)2(\mu_{m}+L^*_4)}ds\nonumber\\
&\leq& 2\int^{\infty}_{t}e^{(s-t)2(\mu_{m}+L^*_4)}L^*_3(
z^+)^2(s)ds-\frac{B_1}{\mu_{m}+L^*_4}.
\end{eqnarray}
Observing (\ref{(bo15)}) and (\ref{(bo16)}), we see that if we prove
one of $( z^+)(t)$ and $( z^-)(t)$ is bounded, the other one can be
deduced to be bounded automatically. Next, we substitute the term $(
z^-)^2(s)$ in (\ref{(bo15)}) by the inequality (\ref{(bo16)}). Then
we can use the change of integration order to get
\begin{eqnarray} \label{(bo17)}
&&( z^+)^2(t)\nonumber\\
&\leq&
2\int^t_{-\infty}e^{(t-s)2(L^*_1-\mu_{m+1})}L^*_2\left[2\int^{\infty}_{s}e^{(r-s)2(\mu_{m}+L^*_4)}L^*_3( z^+)^2(r)dr-\frac{B_1}{\mu_{m}+L^*_4}\right]ds+\frac{A_1}{\mu_{m+1}-L^*_1}\nonumber\\
&\leq&
\lambda\left[\int^t_{-\infty}e^{2(L^*_1-\mu_{m+1})(t-s)}( z^+)^2(s)ds+\int^{\infty}_te^{2(\mu_m+L^*_4)(s-t)}( z^+)^2(s)ds\right]+M,
\end{eqnarray}
where
$$M:=\frac{A_1}{\mu_{m+1}-L^*_1}-\frac{L^*_2B_1}{(\mu_{m+1}-L^*_1)(\mu_m+L^*_4)}>0,$$
and
$$\lambda:=\frac{2L^*_2L^*_3}{\mu_{m+1}-L^*_1-\mu_m-L^*_4}>0.$$
\\
Denote $$\alpha:=\max\{2(\mu_{m+1}-L^*_1),-2(\mu_m+L^*_4)\},\ \ \gamma:=\min\{2(\mu_{m+1}-L^*_1),-2(\mu_m+L^*_4)\}.$$
Then $\alpha, \gamma>0$, and
\begin{eqnarray*}
( z^+)^2(t)&\leq&
M+\lambda\left[\int^t_{-\infty}e^{-\gamma(t-s)}( z^+)^2(s)ds+\int^{\infty}_te^{-\gamma(s-t)}( z^+)^2(s)ds\right].
\end{eqnarray*}
For the above inequality, we consider a variable change for term
$\int^{\infty}_te^{-\gamma(s-t)}( z^+)^2(s)ds$, then
\begin{eqnarray*}
\int^{\infty}_te^{-\gamma(s-t)}( z^+)^2(s)ds
= \int_{-\infty}^{-t}e^{-\gamma(-s-t)}( z^+)^2(-s)ds.
\end{eqnarray*}
Hence
\begin{eqnarray}\label{(bo18)}
(z^+)^2(t)&\leq&
M+\lambda\left[\int^t_{-\infty}e^{-\gamma(t-s)}( z^+)^2(s)ds+\int_{-\infty}^{-t}e^{-\gamma(-s-t)}( z^+)^2(-s)ds\right].
\end{eqnarray}
Replacing $t$ by $-t$ into (\ref{(bo18)}), we have a new form
\begin{eqnarray}\label{(bo19)}
( z^+)^2(-t)&\leq&
M+\lambda\left[\int^{-t}_{-\infty}e^{-\gamma(-t-s)}( z^+)^2(s)ds+\int_{-\infty}^{t}e^{-\gamma(t-s)}( z^+)^2(-s)ds\right].
\end{eqnarray}
Adding (\ref{(bo18)}) and (\ref{(bo19)}) together, we have
\begin{eqnarray*}
( z^+)^2(t)+( z^+)^2(-t)&\leq&
2M+\lambda\bigg[\int^t_{-\infty}e^{-\gamma(t-s)}(( z^+)^2(s)+(
z^+)^2(-s))ds\\
&&\hskip 1.5cm+\int_{-\infty}^{-t}e^{-\gamma(-s-t)}(( z^+)^2(s)+( z^+)^2(-s))ds\bigg].
\end{eqnarray*}
Observing the above inequality, we find that it becomes an induction
problem. Let
$$G(t)=(
z^+)^2(t)+( z^+)^2(-t).$$ Then $G(t)\geq 0$ and
\begin{eqnarray}\label{(bo20)}
G(t)\leq
2M+\lambda\left[\int^t_{-\infty}e^{-\gamma(t-s)}G(s)ds+\int_{-\infty}^{-t}e^{-\gamma(-s-t)}G(s)ds\right].
\end{eqnarray}
To solve this inequality, we use the induction method
by assuming the starting point $G_1(t)\leq 2M$, then
\begin{eqnarray*}
G_1(t)&\leq&2M;\\
G_2(t)&\leq&
2M+\lambda\int^t_{-\infty}e^{-\gamma(t-s)}(2M)ds+\lambda\int_{-\infty}^{-t}e^{-\gamma(-s-t)}(2M)ds\\
&\leq& 2M+2M(\frac{2\lambda}{\gamma});\\
G_3(t)&\leq&
2M+\lambda\int^t_{-\infty}e^{-\gamma(t-s)}(2M+2M(\frac{2\lambda}{\gamma}))ds+\lambda\int_{-\infty}^{-t}e^{-\gamma(-s-t)}(2M+2M(\frac{2\lambda}{\gamma}))ds\\
&\leq& 2M+2M(\frac{2\lambda}{\gamma})+2M(\frac{2\lambda}{\gamma})^2;\\
&\vdots& \\
G_m(t)&\leq&2M+2M(\frac{2\lambda}{\gamma})+2M(\frac{2\lambda}{\gamma})^2+\cdots+2M(\frac{2\lambda}{\gamma})^{m-1};\\
&\vdots&
\end{eqnarray*}
Note
\begin{eqnarray*}
{{2\lambda} \over \gamma}=\frac{\frac{4L^*_2L^*_3}{\mu_{m+1}-L^*_1-\mu_m-L^*_4}}{\gamma}=\frac{8L^*_2L^*_3}{(\alpha+\gamma)\gamma}<1.
\end{eqnarray*}
Hence, $G_m(t)$ has a uniform
bound which does not depend on $N$. This means $( z^+)^2(t)+(
z^+)^2(-t)$ is bounded uniformly in $N$. And since $( z^+)^2(t)$ and
$( z^+)^2(-t)$ must be non-negative, we have $( z^+)^2(t)$ has a
uniform bound. Replacing this bound into (\ref{(bo16)}), we obtain a bound for
$( z^-)^2(t)$. Then we have a bound for $\mid z(t)\mid$. And this bound
completely does not depend on $N$ of $F_N^*$. Hence, we can choose $N$
big enough such that $F_N^*=F^*.$ Then, the globally boundedness
condition for $F$ can be omitted.
$\hfill \hfill $$\sharp$
\medskip

A similar result holds for stochastic differential equation (\ref{(bo50)}) and stochastic integral equation (\ref{(bo1)}).
\begin{thm}
 Assume conditions on $A$, $B_{0}$ in
Theorem \ref{aug20a}. Let $F:  R^d\to  R^d$ be a continuous map satisfying Condition (M) and $\nabla F(u)$ exist and be locally bounded.
Then there exists at least one $\mathcal{F}$-measurable map
$Y:\Omega\rightarrow  R^d$ satisfying (\ref{(bo1)}) and the solution satisfies
\begin{eqnarray*}
u(t, Y(\omega))=Y(\theta_t\omega),\ \ {\rm for \ any} \ \
t\geq 0
\ \ \ \ a.s.
\end{eqnarray*}
\end{thm}

{\bf Acknowledgements}

It is our great pleasure to thank J.Q. Duan, K.D. Elworthy, R. Hudson, K.N. Lu, and J. L. Wu for very useful conversations.   CF would like to acknowledge the support of National Basic Research Program of China
(973 Program No. 2007CB814903), National Natural Science
Foundation of China (No. 70671069 and No. 10971032). BZ wishes to thank the
financial support from Loughborough University.

\end{document}